\documentclass[12pt]{article}
\usepackage{amsmath,amsfonts,amssymb,amsthm,enumerate}

\def\B{\langle}
\def\K{\rangle}

\newcommand{\R}{{\mathbb R}}
\newcommand{\C}{{\mathbb C}}
\newcommand{\Z}{{\mathbb Z}}

\newcommand{\1}{{\bf 1}}
\newcommand{\F}{{\mathbb F}}

\newcommand{\h}{{\bf h}}
\newcommand{\Perm}{}

\newcommand{\al}{\alpha}
\newcommand{\lm}{\lambda}
\newcommand{\Lm}{\Lambda}

\newcommand{\eps}{\varepsilon}

\renewcommand{\th}{\theta}

\newcommand{\BW}{\Lm_{16}}
\newcommand{\ka}{\kappa}

\newcommand{\ch}{{\rm dim_* }}

\newcommand{\Aut}{{\rm Aut}}
\newcommand{\Orbit}{{\mathcal O}}
\newcommand{\AutO}{{O}}
\newcommand{\w}{\omega}

\newcommand{\ot}{\otimes}
\newcommand{\op}{\oplus}
\newcommand{\eqa}{\begin{eqnarray}}
\newcommand{\eeqa}{\end{eqnarray}}
\newcommand{\eqn}{\begin{eqnarray*}}
\newcommand{\eeqn}{\end{eqnarray*}}

\newcommand{\Hom}{{\rm Hom}}
\newcommand{\Irr}{{\rm Irr}}
\newcommand{\Cent}{Z}
\newcommand{\IntType}[3]{\!%
\setbox0=\hbox{$\scriptstyle #3\atop\scriptstyle #1\;#2$}%
{\scriptstyle#3\atop\scriptstyle%
\hbox to\wd0{\hfil$\scriptstyle #1$\hfil\hfil$\scriptstyle#2$\hfil}}%
\!}
\newtheorem{dfn}{Definition}[section]
\newtheorem{pro}[dfn]{Proposition}
\newtheorem{thm}[dfn]{Theorem}

\newtheorem{lem}[dfn]{Lemma}
\newtheorem{cor}[dfn]{Corollary}
\newtheorem{rem}[dfn]{Remark}
\newtheorem{notee}[dfn]{Note}
\newtheorem{exa}[dfn]{Example}

\newcommand{\Imm}{{\rm Im\ }}
\newcommand{\Ker}{{\rm Ker\ }}
\newcommand{\NO}{\,{\raise0.25em\hbox{$\mathop{\hphantom{\cdot}}\limits^{_{\circ}}_{^{\circ}}$}}\,}
\newcommand{\IntTK}[3]{
\scriptstyle\hphantom{\scriptstyle#2}\hfil#1\hfil\hphantom{\scriptstyle#3}%
\atop\scriptstyle#2{\mathop{%
\hbox{$\smash{\rightarrow\joinrel\relbar\joinrel\rightarrow\joinrel\relbar}$}%
}\limits^{\raise5pt\hbox{\kern.3pt${\downarrow}$}\atop\raise4pt\hbox{$\smash{|}$}}}#3%
}
\newcommand{\qe}{\qed\vskip2ex}
\makeatletter
\@addtoreset{equation}{section}

\makeatother

\newcommand{\labell}{\label}
\def\bl{\begin{lem}\labell}
\def\el{\end{lem}}
\def\bt{\begin{thm}\labell}
\def\et{\end{thm}}
\def\bp{\begin{pro}\labell}
\def\ep{\end{pro}}
\def\br{\begin{rem}\labell}
\def\er{\end{rem}}
\def\bc{\begin{cor}\labell}
\def\ec{\end{cor}}
\def\bd{\begin{dfn}\labell}
\def\ed{\end{dfn}}
\def\bn{\begin{notee}\labell}
\def\en{\end{notee}}
\def\be{\begin{exa}\labell}
\def\ee{\end{exa}}
\def\proof{{\it Proof.}}
\def\h{\mathfrak{h}}

\title{
\Large The Automorphism Group\\ of the Vertex Operator Algebra $V_L^+$\\ for an even lattice $L$ without roots}
\author{Hiroki SHIMAKURA\footnote{The author was supported by the Japan Society for the Promotion of Science Research Fellowships for Young Scientists.}}
\date{\small\it Graduate School of Mathematical Sciences,\\
University of Tokyo, Komaba 3-8-1, Tokyo 153-8914, Japan\\
{\rm e-mail: shima@ms.u-tokyo.ac.jp}
}
\begin{document}
\maketitle


\begin{abstract}
The automorphism group of the vertex operator algebra $V_L^+$ is studied by using its action on isomorphism classes of irreducible $V_L^+$-modules.
In particular, the shape of the automorphism group of $V_{L}^+$ is determined when $L$ is isomorphic to an even unimodular lattice without roots, $\sqrt2R$ for an irreducible root lattice $R$ of type $ADE$ and the Barnes-Wall lattice of rank $16$.
\end{abstract}

\section*{Introduction}
Since the introduction of vertex operator algebras (VOAs) \cite{Bo,FLM2}, it has been a basic problem to determine the automorphism group of a VOA.
For example, as is well known, the automorphism group of the moonshine module is the Monster simple group (cf.\ \cite{FLM2}).
We believe that there are many VOAs having an ``interesting" (finite) group as the automorphism group.

The main objective of this paper is the VOA $V_L^+$ which is the fixed-point subspace of the lattice VOA $V_L$ associated with a (positive-definite) even lattice $L$ with respect to an automorphism $\theta_{V_L}$ induced from the $(-1)$ isometry of $L$.
In this paper, we study the automorphism group $\Aut(V_L^+)$ of $V_L^+$.
Clearly $\Aut(V_L^+)$ has the subgroup $H_L$ obtained as the restriction of the centralizer of $\theta_{V_L}$ in $\Aut(V_L)$.
It was proved in \cite{DG} that $\Aut(V_L^+)$ coincides with $H_L$ when $L$ is isomorphic to an even lattice $L$ of rank one except the lattice $2A_1$.
On the other hand, it was shown in \cite{DG,Gr2} that $\Aut(V_L^+)$ does not coincide with $H_L$ when $L$ is isomorphic to $2A_1$ or $\sqrt2E_8$.
We call an automorphism of $V_L^+$ not belonging to $H_L$ an extra automorphism.
So the natural question occurs: ``When does $V_L^+$ have extra automorphisms?".

Our main result is to give a uniform method of determining ${\rm Aut}(V_L^+)$ for an even lattice $L$ without roots (Section \ref{SMAIN}).
As a corollary, we will answer the question above.
More precisely, for an even lattice $L$ without roots, we will show that $V_L^+$ has extra automorphisms if and only if $L$ is obtained from a binary code by Construction B (Proposition \ref{PGOVL+} (2)).
Hence we obtain a new relation between the VOA $V_L^+$ and the lattice $L$.
As an application, we will determine $\Aut(V_L^+)$ when $L$ is isomorphic to an even unimodular lattice without roots, $\sqrt2R$ for an irreducible root lattice $R$ and the Barnes-Wall lattice of rank $16$ (Theorem \ref{autoofunimodular}, \ref{nonsplitofG} and Table \ref{autorootlatticeVOA}).

\medskip

Let us explain our basic idea.
Given a group $G$, one of effective methods of determining the group structure of $G$ is to find a set $X$ on which $G$ acts as permutations and to determine both the kernel and image of the natural group homomorphism from $G$ to the symmetric group $Sym(X)$ on $X$.
In our case, the group $G$ is $\Aut(V_L^+)$, and we take, as such an $X$, the set $S_L$ of all isomorphism classes of irreducible $V_L^+$-modules.

We now recall the action of automorphisms of a VOA on isomorphism classes of its irreducible modules form \cite{DM1}.
Let $(M,Y_M)$ be a module for a VOA $V$.
For an automorphism $g$ of $V$, we set $Y_M^g(v,z)=Y_M(gv,z)$ ($v\in V$).
Then $(M,Y_M^g)$ is a $V$-module.
If $M$ is irreducible then so is $M^g$, hence automorphisms of $V$ acts on the set of all isomorphism classes of irreducible $V$-modules.

We obtained the set $S_L$ on which $\Aut(V_L^+)$ acts.
However, $Sym(S_L)$ is too large to describe $\Aut(V_L^+)$ directly.
So we replace $S_L$ with a subset $P$ of $S_L$ preserved by the action of $\Aut(V_L^+)$, and consider a structure on $P$ associated with the fusion rules.
Since the action of $\Aut(V_L^+)$ preserves the fusion rules, we obtain a group homomorphism from $\Aut(V_L^+)$ to the group consisting of all elements of $Sym(P)$ preserving the structure of $P$.
This homomorphism gives us sufficient information to determine $\Aut(V_L^+)$.

Precisely, we proceed as follows:
\begin{enumerate}[(1)]
\item Describe the stabilizer $H_L$ of the isomorphism class of $V_L^-$ in $\Aut(V_L^+)$.
\item Determine the orbit $Q_L$ of the isomorphism class of $V_L^-$.
\item Find a subset $P$ of $S_L$ such that $P$ contains the orbit $Q_L$ and forms a vector space under the fusion rules.
\end{enumerate}
Here $V_L^-$ is the irreducible $V_L^+$-module which is the $(-1)$-eigenspace of the involution $\theta_{V_L}$.

By performing the steps above, we can determine $\Aut(V_L^+)$ in the following way.
By the step (3), we obtain a natural group homomorphism $\zeta_{V_L^+}$ from $\Aut(V_L^+)$ to the general linear group on the vector space $P$.
So it suffices to determine the kernel and image of $\zeta_{V_L^+}$.
Since the kernel is a subgroup of the stabilizer $H_L$, we can determine it by step (1).
Clearly $\Aut(V_L^+)$ acts transitively on the orbit $Q_L$.
Hence the index of the stabilizer $H_L$ in $\Aut(V_L^+)$ is equal to the cardinality of the orbit $Q_L$, and so is the index of the image of the stabilizer $H_L$ in the image of $\Aut(V_L^+)$ under the homomorphism $\zeta_{V_L^+}$.
By the step (1), we can determine the image of the stabilizer $H_L$.
In this stage, we can use group theoretical facts on general linear groups, and one can compute the image of $\zeta_{V_L^+}$ in principle.

\medskip

Let us explain our arguments in the steps (1), (2) and (3) in detail.
First we describe the stabilizer $H_L$ of the isomorphism class $[0]^-$ of $V_L^-$.
Clearly, any automorphism of $V_L$ that commutes with the involution $\theta_{V_L}$ preserves both $V_L^+$ and $V_L^-$.
Hence by restriction, we obtain a group homomorphism from the centralizer $C_{\Aut(V_L)}(\theta_{V_L})$ of $\theta_{V_L}$ in $\Aut(V_L)$ to $H_L$.
We will show that this homomorphism is surjective, namely we will prove that for any element of $H_L$, there exists its preimage in $C_{\Aut(V_L)}(\theta_{V_L})$ (Theorem \ref{sc} and Proposition \ref{CMain}).
By using the description of $\Aut(V_L)$ given in \cite{DN1}, we can compute the stabilizer $H_L$ of $[0]^-$ in principle and we complete the step (1).

Next we study the orbit $Q_L$ of the isomorphism class $[0]^-$ of $V_L^-$.
When the lattice $L$ does not have roots, then the degree one subspace of $V^+_L$ is zero and the automorphism group of $V^+_L$ is expected to be finite. 
(This is indeed the case as we will show in Proposition \ref{co4}.)
We are interested in such a case, so we assume that $L$ has no roots.
If the isomorphism class $W$ of an irreducible $V_L^+$-module $M$ belongs to $Q_L$, then the following must be satisfied:
\begin{enumerate}[(a)]
\item the dimension of the degree $1$ subspace of $M$ is equal to that of $V_L^-$,
\item the fusion rule $W\times W=V_L^+$ holds.
\end{enumerate}
We note that the isomorphism class of any irreducible $V_L^+$-module is either of untwisted type or of twisted type, which are came from (untwisted) $V_L$-modules or $\theta_{V_L}$-twisted $V_L$-modules respectively (cf. \cite{AD}).
For simplicity, we consider the following two cases:
\begin{enumerate}[(I)]
\item $Q_L$ does not contain isomorphism classes of irreducible $V_L^+$-modules of twisted type,
\item $Q_L$ contains isomorphism classes of irreducible $V_L^+$-modules of twisted type.
\end{enumerate}

In the case (I), let $W$ be the isomorphism class of an irreducible $V_L^+$-module of untwisted type satisfying the conditions $(a)$ and $(b)$, and suppose that $W$ is not the isomorphism class $[0]^-$ of $V_L^-$.
We will prove that the lattice $L$ comes from the Construction B, a well-known procedure which makes a lattice from a binary code, and that $W$ and $[0]^-$ are actually conjugate each other, that is, there exists an element of $\Aut(V_L^+)$ which permutes $W$ and $[0]^-$.
For an even lattice $N$ obtained by Construction B, it was shown in \cite{FLM2} that there exists an automorphism of $V_N^+$ which cannot lift to an automorphism of $V_N$.
In our case, the automorphism of $V_L^+$ permuting $W$ and $[0]^-$ is given by that described as in \cite{FLM2}.
By these consideration, we determine $Q_L$ and complete the step (2) (Theorem \ref{lattice3}).
Then, by using the fusion rules given in \cite{ADL}, we show that the union of $Q_L$ and the isomorphism class of $V_L^+$ forms a vector space over $\F_2$ under the fusion rules (Proposition \ref{lattice5}), and this completes the step (3).

In the case (II), let $W$ be an element of $Q_L$ of twisted type.
By the assumption, the graded dimensions of $V_L^-$ and $W$ are the same, and the fusion rules of $V_L^-$ and those of $W$ are similar.
In this case, by the results in \cite{ADL,FLM2}, we deduce that $L$ is a $2$-elementary totally even lattice of rank $8$ or $16$ (Proposition \ref{orbit3}).
This implies that $L$ is isomorphic to $\sqrt2E_8$ or the Barnes-Wall lattice $\BW$ of rank $16$ under the assumption that $L$ has no roots (cf. \cite{CS,Qu}).
By using the fusion rules of $V_L^+$ given in \cite{ADL}, if $L$ is $2$-elementary totally even then the set $S_L$ of all isomorphism classes of irreducible $V_L^+$-modules forms a vector space over $\F_2$ under the fusion rules.
Furthermore, when the rank of $L$ is a multiple of $8$, we define a natural $\Aut(V_L^+)$-invariant quadratic form on this vector space (Theorem \ref{PQuad}), and obtain a group homomorphism from $\Aut(V_L^+)$ to an orthogonal group.
By using this homomorphism, we determine $\Aut(V_L^+)$ when $L$ is isomorphic to $\sqrt2E_8$ and $\BW$.
As a corollary, we see that these lattices satisfy the condition (II).

\medskip

The organization of this paper is as follows:
In Section \ref{ChVOA}, we recall some definitions and facts necessary in this paper.
In Section \ref{ChVLM}, we review the properties of the VOA $V_L^+$, such as the irreducible modules, the fusion rules, the graded dimensions and automorphisms.
In Section \ref{ChAGVL+}, we study the action of $\Aut(V_L^+)$ on the set of all isomorphism classes of irreducible $V_L^+$-modules, and give a method of determining the automorphism group of $V_L^+$ for an even lattice $L$ without roots.
In Section \ref{SExample}, applying this method to many important lattices $L$, we determine $\Aut(V_L^+)$.

\medskip

Throughout this paper, we will work over the field $\C$ of complex numbers unless otherwise stated.
We denote by $\Z_q$ the set of integers by $\Z$ and the ring of integers modulo $q$ .
We often identify $\Z_2$ with the field $\F_2$ of two elements.
We denote the elementary abelian $2$-group of order $2^n$ simply by $2^n$ and the root lattice of type $R$ simply by $R$.
We use the notation of \cite{ATLAS} for classical groups except that we denote by $O_n^+(2)$ the orthogonal group of degree $n$ over $\F_2$ of plus type and by $\Omega_n^+(2)$ its simple normal subgroup of index $2$ .
We also use the notation of \cite{ATLAS} for group extensions:
For groups $A$ and $B$, $A.B$ denotes any group extension for which the quotient by a normal subgroup isomorphic to $A$ is isomorphic to $B$, $A:B$ indicates any case of $A.B$ for which the extension is split and $A\cdot B$ indicates any case of $A.B$ for which the extension is not split.

\medskip
 
{\it Acknowledgments.} The author would like to thank Professor Atsushi Matsuo for giving helpful advice and comments on this paper.
He also thanks Professor Alexander A. Ivanov for stimulating discussions on orthogonal groups, Professor Satoshi Yoshiara for valuable comments on finite groups, Hiroshi Yamauchi for useful advice on the simple current extension of a simple VOA and Toshiyuki Abe for sending him a preprint concerning the subject.
Finally he thanks the referee for the valuable suggestions.
\medskip

\section{Preliminaries}\label{ChVOA}
In this section, we recall or give some definitions and facts necessary in this paper.

\subsection{Fusion rules, graded dimensions and automorphisms}\label{SATA}
In this subsection, we recall the fusion rules and the graded dimensions from \cite{FLM2,FHL} and the action of automorphisms of a VOA on its modules from \cite{DM1}.
For details of the definition of VOAs, see \cite{Bo,FLM2}.

Let $V$ be a VOA and let $M^1$, $M^2$, $M^3$ be $V$-modules.
Let $I_{V}\left(\IntType{M^1}{M^2}{M^3}\right)$ denote the vector space generated by all intertwining operators of type $\left(\IntType{M^1}{M^2}{M^3}\right)$ (for the definition see \cite{FHL}).
We set 
\begin{eqnarray*}
N_{M^{1}M^{2}}^{M^{3}}=\dim I_{V}\left(\IntType{M^1}{M^2}{M^3}\right).
\end{eqnarray*}
This number is called the {\it fusion rule}.
The following symmetry on the fusion rules is well known.
\bp{PFHL}{\rm \cite{FHL}} Let $M^1,M^2,M^3$ be $V$-modules.
Then
\begin{eqnarray*}
N_{M^{1}M^{2}}^{M^{3}}=N_{M^{2}M^{1}}^{M^{3}}.
\end{eqnarray*}
\ep

Let $\{W^i\}_{i\in I}$ be the set of all isomorphism classes of irreducible $V$-modules and let $M^i$ be a representative of $W^i$ for $i\in I$.
We denote the fusion rules by using the following formal product:
\begin{eqnarray}
W^i\times W^j=\sum_{k\in I}N_{M^{i}M^{j}}^{M^{k}}W^k.\label{FormalP}
\end{eqnarray}
We note that this product is independent of the choice of representatives, and that $W^i\times W^j=W^j\times W^i$ by Proposition \ref{PFHL}.
Later, when $V$ is the VOA $V_L^+$, we will find a subset $I^\prime$ of $I$ such that for any $i,j\in I^\prime$ there exists a unique element $k\in I^\prime$ satisfying $W^i\times W^j=W^k$, and will view this formal product as a binary operation on $\{W^i\}_{i\in I^\prime}$.
We often denote the isomorphism class of an irreducible $V$-module $M$ simply by $M$.

Let $B=\oplus_{i\in\C}B_{i}$ be a graded vector space over $\C$ satisfying $\dim B_i<\infty$ for all $i\in\C$.
Then the {\it graded dimension} of $B$ is the formal sum
\begin{eqnarray*}
\ch B=\sum_{i\in\C}(\dim B_i)q^i.
\end{eqnarray*}
Since any module of a VOA $V$ is a graded vector space whose homogeneous subspaces are finite dimensional, we may consider its graded dimension.
Clearly, if $V$-modules $M^1$ and $M^2$ are isomorphic then $\ch M^1=\ch M^2$.
Hence for the isomorphism class $W$ of a $V$-module $M$, we often denote by $\ch W$ the graded dimension of $M$.

\bn{NCH}{\rm The graded dimension is often called the character.
However, since we consider characters of groups later, we use the term graded dimension in order to avoid confusion.}
\en

An {\it automorphism} of a VOA $(V,Y,\1,\w)$ is a linear isomorphism $g:V\to V$ satisfying $Y(gv,z)g=gY(v,z)$ for all $v\in V$ that fixes the Virasoro element $\w$.
Let $\Aut(V)$ denote the group of all automorphisms of $V$.
It is easy to see that any automorphism preserves the grading and fixes the vacuum vector $\1$.

Now we recall the action of automorphisms of $V$ on its modules.
Let $(M,Y_M)$ be a $V$-module.
For an automorphism $g$ of $V$, the linear map $Y_{M^g}$ is defined by
\begin{eqnarray*}
Y_{M^g}(v,z)=Y_M(gv,z).
\end{eqnarray*}
Then $M^g=(M,Y_{M^g})$ is a $V$-module.
The following proposition concerning $M^g$ is an easy exercise.
\bp{RAM}
\begin{enumerate}
\item If $M$ is irreducible then $M^g$ is also irreducible.
In particular $\Aut(V)$ acts on the set of all isomorphism classes of irreducible $V$-modules.
\item The graded dimensions of $M$ and $M^g$ are the same.
\item Suppose $V_0=\C\1$.
Then, for $v\in V_1$, $g=\exp(v_0)$ is an automorphism of $V$ and $M^g\cong M$ as $V$-modules.
\item Let $M^1$, $M^2$, $M^3$ be irreducible $V$-modules.
Then $N_{M^1,M^2}^{M^3}=N_{M^{1,g},M^{2,g}}^{M^{3,g}}$.
\end{enumerate}
\ep

\subsection{Construction B and $2$-elementary totally even lattices}\label{SLC}
In this subsection, we study lattices obtained by Construction B from binary codes and $2$-elementary totally even lattices.
For details of the definitions of lattices and codes, see \cite{CS}.

Let $n$ be a positive integer and set $\Omega_n=\{1,2,\dots,n\}$.
Let $\{\alpha_i|\ i\in\Omega_n\}$ be an orthogonal basis of $\R^n$ satisfying $\langle\alpha_i,\alpha_j\rangle=2\delta_{i,j}$.
For a subset $J\subset\Omega_n$, we set $\alpha_J=\sum_{i\in J}\alpha_i$.
Let $C$ be a (binary linear) code of length $n$, which we regard as a subset of the power set of $\Omega_n$.
Then the lattice
\begin{eqnarray*}
L_B(C)=\sum_{c\in C}\Z\frac{1}{2}\alpha_c+\sum_{i,j\in\Omega_n}\Z(\alpha_i+\alpha_j)
\end{eqnarray*}
is called {\it the lattice obtained by Construction B from $C$} associated with $\{\alpha_i|\ i\in\Omega_n\}$.
The following propositions concerning $L_B(C)$ are well-known.

\bp{PGCB}Let $\tilde{C}$ be a set of generators of $C$.
Then $L_B(C)$ is generated by $\{\alpha_i\pm \alpha_j,\ \frac{1}{2}\alpha_c|\ i,j\in\Omega_n,\ c\in \tilde{C}\}$ as a $\Z$-module.
\ep

\bp{PCOB}{\rm \cite{CS}} Let $C$ be a doubly even code of length $n$.
\begin{enumerate}
\item The lattice $L_B(C)$ is even and its determinant is $2^{n-2\dim C+2}$.
\item The dual lattice $L_B(C)^\circ$ of $L_B(C)$ is given by
\begin{eqnarray*}
L_B(C)^\circ=L_B(C^\perp)+\Z\alpha_1+\Z\frac{\alpha_{\Omega_n}}{4}.
\end{eqnarray*}
\item If the minimum weight of $C$ is greater than $4$ then the minimum norm of $L_B(C)$ is $4$, and in particular $L_B(C)$ has no roots.
\end{enumerate}
\ep

The following are typical examples of lattices obtained by Construction B.

\be{REXCB}{\rm \cite{CS}
\begin{enumerate}\item Let $C$ be the code of length $n$ consisting of the all-zero codeword.
Then $L_B(C)$ is isomorphic to $\sqrt2D_n$, where $D_1$, $D_2$ and $D_3$ are regarded as $2A_1$, $A_1\oplus A_1$ and $A_3$ respectively.
\item Let $C$ be the code of length $8$ consisting of the all-zero and all-one codewords.
Then $L_B(C)$ is isomorphic to $\sqrt2E_8$.
\item Let $RM(4,1)$ denote the first order Reed-Muller code of length $2^4$.
Then $L_B(RM(4,1))$ is isomorphic to the Barnes-Wall lattice of rank $16$.
\end{enumerate}}
\ee

\br{RARLBW}{\rm The automorphism groups and the determinants of the irreducible root lattices and the Barnes-Wall lattice $\BW$ of rank $16$ are summarized in Table \ref{TaAR} (cf. Chapter 4 of \cite{CS}).}
\er

\begin{table}
\caption{Automorphism groups of lattices}\label{TaAR}
\begin{center}
\begin{tabular}{|c|| c|c|} 
\hline 
$R$&$\det R$&$\Aut(R)$\\ \hline
$A_1$&$2$&$\Z_2$\\ \hline
$A_n$ $(n\neq1)$&$n+1$&${\Perm Sym}_{n+1}\times 2$\\ \hline
$D_4$&$4$&$(2^{3}:{\Perm Sym}_4):{\Perm Sym}_3$\\ \hline
$D_n$ $(n\ge5)$&$4$&$2^{n-1}:{\Perm Sym}_n:2$\\ \hline
$E_6$&$3$&$2.U_4(2):2$\\ \hline
$E_7$&$2$&$2.Sp_6(2)$\\ \hline
$E_8$&$1$&$2.O_{8}^+(2)$\\ \hline
$\BW$&$2^8$&$2_+^{1+8}\cdot\Omega_8^+(2)$\\
\hline
\end{tabular}
\end{center}
\end{table}

The following proposition will be used in Section \ref{SAGVL+}.
For a subset $U$ of $\R^n$ and $i\in\Z_{>0}$, we denote $U_i=\{v\in U|\ \langle v,v\rangle=i\}$.

\bp{lattice2} Let $L$ be an even lattice of rank $n$ without roots.
If there exists a vector $\lambda\in L^\circ\cap (L/2)$ such that $|(\lambda+L)_2|=2n$ then $L$ is obtained by Construction B associated with an orthogonal basis of $\R\otimes_{\Z}L$ consisting of vectors in $(\lambda+L)_2$.
\ep

\proof\ First, let us show that $(\lambda+L)_2$ forms a root system of type $A_1^n$.
Since $\langle\lambda,\lambda\rangle\in2\Z$ and $\lambda\in L^\circ$, the lattice $N=L+\Z\lambda=L\cup(\lambda+L)$ is even.
It follows from $L_2=\phi$ that $(\lambda+L)_2$ forms a root system.
Assume that it is not of type $A_1^n$.
Since $|(\lambda+L)_2|$ is equal to the number of roots of $A_1^n$, the set $(\lambda+L)_2$ contains a root system not of type $A_1$.
Hence there exist $x,y\in(\lambda+L)_2$ such that $\langle x,y\rangle=-1$.
Then the norm of $x+y$ is $2$ and $x+y\in L$, which is a contradiction for $L_2=\phi$.
Thus $(\lambda+L)_2$ forms a root system of type $A_1^n$.
In particular $(\lambda+L)_2$ contains an orthogonal basis of $\R\otimes_{\Z}L$.

Next we show that $L$ is obtained by Construction B.
Let $\{\alpha_i|\ i\in\Omega_n\}$ be a basis of $\R\otimes L$ consisting of vectors in $(\lambda+L)_2$ and set $M=\oplus_{i\in\Omega_n}\Z\alpha_i$.
Then $M^\circ=M/2$ and $M\subset N\subset M^\circ$.
Let $C$ be the code obtained as the image of $N/M$ under the homomorphism of vector spaces over $\F_2$: $M^\circ/M\to\Z_2^n$, $\sum_{i\in\Omega_n}c_i\alpha_i/2\mapsto (\bar{c}_1,\bar{c}_2,\dots,\bar{c}_n)$.
Then we have $N=\sum_{i\in\Omega_n}\Z\alpha_i+\sum_{c\in C}\Z\alpha_c/2$.
It follows from $\{\alpha_i|\ i\in\Omega_n\}\subset\lambda+L$ and $2\lambda\in L$ that $\alpha_i\pm\alpha_j\in L$.
Let $\{c_1,c_2,\dots,c_k\}$ be a basis of $C$.
Then there exists $D\subset\Omega_n$ such that $\alpha_{c_i}/2$ belongs to $\epsilon_D(L)$ for all $i$, where $\epsilon_D$ denotes the linear map of $\R^n$ given by $\epsilon_D(\alpha_i)=\alpha_i$ if $i\in D$, and $\epsilon_D(\alpha_i)=-\alpha_i$ if $i\notin D$.
Since $L$ contains $\alpha_i\pm\alpha_j$ for all $i,j\in\Omega_n$, so does $\epsilon_D(L)$.
Hence it follows from Proposition \ref{PGCB} that $\epsilon_D(L)\supset L_B(C)$.
Since $|N/L_B(C)|=|N/L|=2$, we have $\epsilon_D(L)= L_B(C)$.
Thus $L$ is obtained by Construction B from $C$ associated with $\{\epsilon_D(\alpha_i)|\ i\in\Omega_n\}$.
\qe

In the rest of this section, we consider $2$-elementary totally even lattices.
The VOA $V_L^+$ for such a lattice will be studied in Section \ref{SVL+2ETE}

An even lattice $L$ is called {\it $2$-elementary totally even} if $2L^\circ\subset L$ and if $\sqrt2L^\circ$ is even.
Then the lattice $\sqrt2E_8$ and the Barnes-Wall lattice of rank $16$ are characterized as follows:
\bp{PEX2ETE} 
\begin{enumerate}
\item {\rm (cf. \cite{CS})} The lattice $\sqrt2E_8$ is the unique $2$-elementary totally even lattice of rank $8$ with determinant $2^8$ up to isomorphism.
\item {\rm \cite[Theorem 4]{Qu}} The Barnes-Wall lattice of rank $16$ is {the} unique $2$-elementary totally even lattice of rank $16$ with determinant $2^8$ without roots up to isomorphism.
\end{enumerate}
\ep

\section{Vertex operator algebra $V_L^+$}\label{ChVLM}
In this section, we review the VOA $V_L^+$ and its automorphisms.

\subsection{$V_L^+$ and its irreducible modules}\label{SLVOA}
In this subsection, we review the VOA $V_L^+$ and its irreducible modules as well as their properties from \cite{FLM2,DN2,AD,Ab,ADL}.

Let $L$ be an even lattice of rank $n$ equipped with a positive-definite symmetric $\Z$-bilinear form $\B\cdot,\cdot\K$.
We set $\h_L=\C\otimes_{\Z} L$ and regard $L$ as a subgroup of $\h_L$.
We extend the form $\B\cdot,\cdot\K$ to a $\C$-bilinear form on $\h_L$.
Let $L^\circ=\{v\in\h_L|\ \B v,L\K\subset\Z\}$ denote the dual lattice of $L$.
For a subset $M\subset L^\circ$, we set $\h_M=\C\B M\K\subset\h_L$.

Let $q$ be the minimal positive even integer such that $\B L^\circ,L^\circ\K\subset2\Z/q$.
By Remark 12.18 in \cite{DL}, there exists a central extension $\hat{L^\circ}$ of $L^\circ$ by the cyclic group $\B\omega_q\K\cong\Z_q$
\begin{eqnarray} 1\rightarrow\B\omega_q\rangle\rightarrow\hat{L^\circ}~\bar{\to}~L^\circ\rightarrow 0\label{cent}
\end{eqnarray}
such that $[a,b]=\kappa_L^{\B \bar{a},\bar{b}\K}$ if $\bar{a},\bar{b}\in L$, where $\kappa_L=\w_q^{q/2}$.
We often regard $\omega_q$ as a primitive $q$-th root of the unity.
We may choose a $\Z$-bilinear $2$-cocycle $\eps(\cdot,\cdot):L^\circ\times L^\circ\to\Z_q$ such that $\eps(\alpha,\alpha)=q\langle\alpha,\alpha\rangle/4$ and denote the corresponding section by $L^\circ\to \hat{L^\circ}$, $\alpha\mapsto e(\alpha)$.
Let $\hat{L}$ be a central extension of $L$ by $\langle\kappa_L\rangle\cong\Z_2$ such that $[a,b]=\kappa_L^{\B \bar{a},\bar{b}\K}$ for $a,b\in\hat{L}$.
We regard $\hat{L}$ as a subgroup of $\hat{L^\circ}$.

Let us consider the lattice VOA $V_L$.
Form the induced $\hat{L^\circ}$-module $\C\{L^\circ\}=\C[\hat{L^\circ}]\ot_{\C[\omega_q]}\C$, where $\C[\cdot]$ denotes the operation of forming the group algebra and $\omega_q$ acts on $\C$ by the multiplication.
For a subset $M\subset L^\circ$, we set $\hat{M}=\{a\in \hat{L^\circ}|\ \bar{a}\in M\}$, $\C\{M\}=\C[\hat{M}]\ot_{\C[\omega_q]}\C$ and $V_{M}=S(\h_M\ot t^{-1}\C[t^{-1}])\ot\C\{M\}$.
For convenience, we denote $h\ot t^m$ by $h(m)$ for $h\in \h_M$ and $m\in\Z$.
It was shown in \cite{Bo,FLM2} that $V_L$ has a simple VOA structure and $V_{\lm+L}$ has an irreducible $V_L$-module structure for $\lm+L\in L^\circ/L$.

The VOA $V_L^+$ is constructed as follows.
Let $\th_{L^\circ}$ be an automorphism of $L^\circ$ such that $\th_{L^\circ}(\omega_q^me(\alpha))=\omega_q^me(-\alpha)$, and let $\th_{V_L}$ be the unique automorphism of the $V_L$-module $V_{L^\circ}$ such that $\th_{V_L}(h(m)\ot1)=-h(m)\ot1$ and $\th_{V_L}(1\ot e(\alpha))=1\ot \th_{L^\circ}(e(\alpha))$ for $h\in \h_{L^\circ}$ and $\al\in L^\circ$.
For a $\th_{V_L}$-stable subspace $E$, we set $E^\pm=\{v\in E|\ \th_{V_L}(v)=\pm v\}$.
Then $V_L^+$ is a subVOA of $V_L$, and $V_{\mu+L}$ and $V_{\lm+L}^\pm$ are irreducible $V_L^+$-modules for $\mu\in L^\circ\setminus (L/2)$ and $\lm\in L^\circ\cap(L/2)$.
Such an irreducible $V_L^+$-module is said to be {\it of untwisted type}

Let us consider irreducible $V_L^+$-modules of twisted type.
Set $K_{L}=\{a^{-1}\th_{L^\circ}(a)|\ a\in \hat{L}\}$.
Then $K_L$ is a normal subgroup of $\hat{L}$.
For an $\hat{L}/K_L$-module $T$ on which $\ka K_L$ acts by $-1$, we set $V_L^{T}=S(\h_L\ot t^{-1/2}\C[t^{-1}])\ot T$.
It was shown in \cite{FLM2} that $V_L^T$ has a $V_L^+$-module structure.
Let $\phi$ be the unique commutative algebra automorphism of $S(\h_L\ot t^{-1/2}\C[t^{-1}])$ such that $\phi(h(m))=-h(m)$ for $h\in \h_L$ and $m\in1/2+\Z$, and let $\th_{V_L^T}$ be the automorphism of the $V_L^+$-module $V_L^T$ given by $\th_{V_L^T}(x\ot t)=\phi(x)\ot t$.
We also denote the $\pm1$-eigenspace of $\th_{V_L^T}$ by $V_L^{T,\pm}$.
Then $V_L^{T,\pm}$ are irreducible $V_L^+$-modules if $T$ is irreducible.
Such an irreducible $V_L^+$-module is said to be {\it of twisted type}.

We study irreducible $\hat{L}/K_L$-modules more precisely.
For a subgroup $G$ of $\hat{L}/K_L$ containing $\langle\kappa_L K_L\rangle$, we denote
\begin{eqnarray*}
X(G)=\{\chi\in\Irr(G)|\ \chi(\kappa_L K_L)=-1\},
\end{eqnarray*}
where $\Irr(G)$ is the set of all irreducible characters of the group $G$.
Applying Theorem 5.5.1 of \cite{FLM2} to $\hat{L}/K_L$, we obtain the following proposition.

\bp{PTC} {\rm \cite{FLM2}} The restriction map $X(\hat{L}/K_L)\to X(Z(\hat{L}/K_L))$ is bijective, where $Z(\hat{L}/K_L)$ is the center of $\hat{L}/K_L$.
Moreover the degree of any element of $X(\hat{L}/K_L)$ is equal to $|L/(L\cap2L^\circ)|^{1/2}$.
\ep 
We denote the irreducible $\hat{L}/K_L$-module corresponding to $\chi\in X(\Cent(\hat{L}/K_L))$ by $T_{\chi}$.

\medskip

We have obtained some irreducible $V_L^+$-modules.
In \cite{DN2} (for rank one) and \cite{AD} (for general rank), the irreducible $V_L^+$-modules were classified.
\bt{Eq:CVL+} {\rm \cite{DN2,AD}} Let $L$ be an even lattice.
Then any irreducible $V_L^+$-module is isomorphic to one of $V_{\lambda+L}^\pm$ $(\lambda\in L^\circ\cap (L/2))$, $V_{\mu+L}$ $(\mu\in L^\circ\setminus (L/2))$ and $V_L^{T_{\chi},\pm}$ $(\chi\in X(\Cent(\hat{L}/K_L)))$.
\et

\br{Risomod}{\rm The irreducible $V_L^+$-modules in Theorem \ref{Eq:CVL+} are non-isomorphic except that $V_{\mu+L}\cong V_{-\mu+L}$ for $\mu\in L^\circ\setminus (L/2)$.}
\er

In this paper, we use the following notation of \cite{ADL} for the isomorphism classes of irreducible $V_L^+$-modules: $[\lambda]^\pm$ $(\lambda\in L^\circ\cap (L/2))$, $[\mu]$ $(\mu\in L^\circ\setminus (L/2))$ and $[\chi]^\pm$ $(\chi\in X(\Cent(\hat{L}/K_L)))$ denote the isomorphism classes of $V_{\lambda+L}^\pm$, $V_{\mu+L}$ and $V_{L}^{T_\chi,\pm}$ respectively

In \cite{Ab} (for rank one) and \cite{ADL} (for general rank), the fusion rules of $V_L^+$ were completely determined.
In particular, we have the following proposition.

\bp{RFusionFo}{\rm \cite{Ab,ADL}}\begin{enumerate}
\item Let $\lambda_1,\lambda_2\in L^\circ\cap(L/2)$ and $\mu\in L^\circ\setminus (L/2)$.
Let $W$ be the isomorphism class of an irreducible $V_L^+$-module.
Then 
\begin{eqnarray*}
{}[0]^-\times & [\lambda_1]^\pm &=[\lambda_1]^\mp,\\ 
{}[0]^-\times & [\mu] &=[\mu],\\
{}[\lambda_1]^{+}\times & [\lambda_2]^{+} &=[\lambda_1+\lambda_2]^{\varepsilon}\ \text{\rm for some}\ \varepsilon\in\{\pm\},\\
{}[0]^+\times & W &=W,\\
{}[\mu]\times & [\mu] &=[0]^++[0]^-+[2\mu].
\end{eqnarray*}
\item Let $M^1,M^2,M^3$ be irreducible $V_L^+$-modules.
Suppose that $M^1$ is of twisted type and that $N_{M^1,M^2}^{M^3}=1$.
If $M^2$ is of twisted type (resp.\ of untwisted type) then $M^3$ is of untwisted type (resp.\ of twisted type).
\item The formal product $\times$ given in (\ref{FormalP}) is associative.
\end{enumerate}
\ep

The following proposition summarizes the graded dimensions of irreducible $V_L^+$-modules (cf.\ Section 10.5 of \cite{FLM2}).
\bp{PCH} Let $L$ be an even lattice of rank $n$.
Let $\lambda\in L^\circ\cap (L/2)$, $\mu\in L^\circ\setminus (L/2)$ and $\chi\in X(Z(\hat{L}/K_L))$ such that $\lambda\notin L$.
Then we have
\begin{eqnarray*}
\ch[\mu]&=&\frac{\Theta_{\mu+L}(q)}{\Phi(q)^n},\\
\ch[0]^\pm&=&\frac{1}{2}\Big(\frac{\Theta_L(q)}{\Phi(q)^n}\pm\frac{\Phi(q)^n}{\Phi(q^2)^n}\Big),\\
\ch[\lambda]^\pm&=&\frac{\Theta_{\lambda+L}(q)}{2\Phi(q)^n},\\
\ch[\chi]^\pm&=&\frac{\dim T_\chi q^{\frac{n}{16}}}{2}\Big(\frac{\Phi(q)^n}{\Phi(q^{1/2})^n}\pm\frac{\Phi(q^2)^n\Phi(q^{1/2})^n}{\Phi(q)^{2n}}\Big),
\end{eqnarray*}
where $\Phi(q)=\Pi_{i=1}^\infty(1-q^i)$ and $\Theta_{U}(q)=\sum_{\alpha\in U}q^{\langle \alpha,\alpha\rangle/2}$ for $U\subset\R^n$.
In particular, first few coefficients are given as follows:
\begin{eqnarray*}
\ch[0]^+&\in&1+\frac{|L_2|}{2}q+\Big(\frac{n|L_2|+|L_4|}{2}+{n+1\choose2}\Big)q^2+q^3\Z[q],\\
\ch[0]^-&\in&\Big(\frac{|L_2|}{2}+n\Big)q+\Big(\frac{n|L_2|+|L_4|}{2}+n\Big)q^2+q^3\Z[q],\\
\ch[\lambda]^\pm&\in&\frac{|(\lambda+L)_{\iota(\lambda+L)}|}{2}q^{\iota(\lambda+L)/2}+q^{\iota(\lambda+L)/2+1}\Z[q],\\
\ch[\chi]^+&\in&\dim T_\chi q^{\frac{n}{16}}+{n+1\choose2}\dim T_\chi q^{\frac{n}{16}+1}+q^{\frac{n}{16}+2}\Z[q],\\
\ch[\chi]^-&\in&n\dim T_\chi q^{\frac{n+8}{16}}+q^{\frac{n+8}{16}+1}\Z[q],
\end{eqnarray*}
where $\iota(\lambda+L)=\min\{\langle v,v\rangle|\ v\in\lambda+L\}$.
\ep

\subsection{Automorphism group of $V_L$}\label{SLAL}
In this subsection, we describe the automorphism group of $V_L$ after \cite{DN1}.
In particular, we study the action of some automorphisms of $V_L$ preserving $V_L^+$ on the set of isomorphism classes of irreducible $V_L^+$-modules.

Let $L$ be an even lattice of rank $n$.
Let $\Aut(\hat{L^\circ})$ denote the automorphism group of the group $\hat{L^\circ}$.
For $g\in\Aut(\hat{L^\circ})$, let $\bar{g}$ denote the linear isomorphism of $L^\circ$ defined by $\bar{g}(\alpha)=\overline{g(e(\alpha))}$, $\alpha\in L^\circ$.
Let $\AutO(L^\circ)$ denote the group of automorphisms of $L^\circ$ preserving the inner product $\langle\cdot,\cdot\rangle$.
Set $\AutO(\hat{L^\circ})=\{g\in\Aut(\hat{L^\circ})|\ g\omega_q=\omega_q,\ \bar{g}\in\AutO(L^\circ)\}$.
We view an element $f\in\Hom(L^\circ,\Z_q)$ as the automorphism of $\hat{L^\circ}$ which sends $e(\alpha)$ to $\omega_q^{f(\alpha)}e(\alpha)$.
Hence we obtain an embedding $\Hom(L,\Z_q)\subset \AutO(\hat{L^\circ})$.
Similarly, the group $O(\hat{L})$, a group homomorphism from $O(\hat{L})$ to $O(L)$ and an embedding $\Hom(L,\Z_2)\subset \AutO(\hat{L})$ are defined.
By Proposition 5.4.1 of \cite{FLM2}, we have the following exact sequences:
\begin{eqnarray}
1\rightarrow\Hom(L^\circ,\Z_q)\hookrightarrow \AutO(\hat{L^\circ})~\bar{\rightarrow}~\AutO(L^\circ)\rightarrow 1,\label{exactauto}\\
1\rightarrow\Hom(L,\Z_2)\hookrightarrow \AutO(\hat{L})~\bar{\rightarrow}~\AutO(L)\rightarrow 1.\label{exactauto2}
\end{eqnarray}

\br{Rsplit}{\rm If $\langle L,L\rangle\subset2\Z$ then the exact sequence (\ref{exactauto2}) is split.}
\er

For $\beta\in L^\circ/2L^\circ$, let $f_\beta:L\to\Z_2$ denote the element of $\Hom(L,\Z_2)$ given by
\begin{eqnarray}
f_\beta:\gamma\mapsto \langle \beta,\gamma\rangle\mod2.\label{Def:fbeta}
\end{eqnarray}
It is easy to see that $\Hom(L,\Z_2)=\{f_\beta|\ \beta\in L^\circ/2L^\circ\}$.

Since $O(\hat{L})$ fixes $\kappa_L$, an element of $O(\hat{L})$ is extended to an automorphism of $\C\{L\}$.
Then $g\in\AutO(\hat{L})$ is also extended to a linear automorphism of $V_L$ by setting
\begin{eqnarray}
g(h_1(n_1)\cdots h_k(n_k)\otimes a)=\bar{g}(h_1)(n_1)\cdots\bar{g}(h_k)(n_k)\otimes g(a),\label{act}
\end{eqnarray}
where $n_i\in\Z_{<0}$, $h_i\in\h_L$ and $a\in\C\{L\}$.
By \cite[Corollary 10.4.8]{FLM2}, this extension gives an automorphism of the VOA $V_L$.
Thus we have an injective group homomorphism $O(\hat{L})\hookrightarrow \Aut(V_L)$ of which the image is identified with $O(\hat{L})$.
For $g\in \AutO(\hat{L})$, we call $g$ a {\it lift of $\bar{g}$}.
By \cite{DN1}, the automorphism group  $\Aut(V_L)$ of $V_L$ is described as follows:
\bt{PAVL}{\rm \cite[Theorem 2.1]{DN1}} Let $L$ be an even lattice.
Then $\Aut(V_L)=N(V_L)O(\hat{L})$, where $N(V_L)=\langle\exp(v_0)|\ v\in (V_L)_1\rangle$ is a normal subgroup of $\Aut(V_L)$.
Moreover, $\Aut(V_L)/N(V_L)$ is isomorphic to a quotient group of $O(L)$.
\et

In the previous section, the involution $\theta_{V_L}$ of $V_L$ is given to define the VOA $V_L^+$.
Clearly $\theta_{V_L}$ is a lift of the $(-1)$-isometry which sends $\alpha\in L$ to $-\alpha$.
The centralizer of $\theta_{V_L}$ in $\Aut(V_L)$ plays an important role in this paper.

\bp{co1}The centralizer $C_{\Aut(V_L)}(\theta_{V_L})$ of $\theta_{V_L}$ in $\Aut(V_L)$ is $C_{N(V_L)}(\theta_{V_L})O(\hat{L})$.
In particular, if $L_2=\phi$ then $C_{\Aut(V_L)}(\theta_{V_L})=\AutO(\hat{L})\cong (\Z_2)^n.O(L)$.
\ep
\proof\ By Theorem \ref{PAVL}, $\Aut(V_L)=N(V_L)O(\hat{L})$.
Since the $(-1)$-isometry belongs to $Z(O(L))$ and $\theta_{V_L}$ commutes with $\Hom(L,\Z_2)$, we have $O(\hat{L})\subset C_{\Aut(V_L)}(\theta_{V_L})$.
Since $N(V_L)$ is normal, the first assertion follows.

Suppose $L_2=\phi$.
Then we have $(V_L)_1=\h_L(-1)$.
It suffices to show that $C_{N(V_L)}(\theta_{V_L})\subset O(\hat{L})$.
Let $h(-1)\in\h_L(-1)$.
Then the exponential $\exp({h(-1)_0})$ sends $x\otimes a\in V_L$ to $\exp({\langle h,\bar{a}\rangle})x\ot a$.
Hence $\exp({h(-1)_0})$ commutes with $\theta_{V_L}$ if and only if $h\in \pi \sqrt{-1}L^\circ$.
It is easy to see that $\exp({\pi\sqrt{-1}\beta(-1)_0})=f_\beta$ for $\beta\in L^\circ$.
Thus $C_{N(V_L)}(\theta_{V_L})=\Hom(L,\Z_2)\subset O(\hat{L})$.
\qe

By the definition of $V_L^+$ and the proposition above, the group $C_{\Aut(V_L)}(\theta_{V_L})$ acts on $V_L^+$.
We describe the action of $O(\hat{L})$ on the set of all isomorphism classes of irreducible $V_L^+$-modules of untwisted type.

\bp{actioncoset} For $g\in\AutO(\hat{L})$, we have 
\eqn
[\mu]^g&=&[\bar{g}(\mu)],\ \mu\in L^\circ\setminus (L/2),\\ 
{}\{[\lambda]^{\pm,g}\}&=&\{[\bar{g}(\lambda)]^\pm\},\ \lambda\in L^\circ\cap(L/2),\\
{}[0]^{\pm,g}&=&[0]^\pm,
\eeqn
where we regard $\bar{g}\in \AutO(L)$ as an automorphism of $L^\circ$.
Moreover, if $\bar{g}=1$ then for $\lambda\in L^\circ\cap(L/2)$,
\eqn
[\lambda]^{\pm,g}=
\left\{\begin{array}{cl}
 \mbox{$[\lambda]^\pm$} & \mbox{${\rm if}\ g(2\lambda)=0$},\\
 \mbox{$[\lambda]^\mp$} & \mbox{${\rm if}\ g(2\lambda)=1$},
\end{array}
\right.
\eeqn
where we regard $g$ as an element of $\Hom(L,\Z_2)$.
\ep
\proof\ Let $\alpha$ be a vector in $L^\circ$.
Since $\AutO(L)=\AutO(L^\circ)$, there exists an element $h$ of $\AutO(\hat{L^\circ})$ such that $\bar{g}=\bar{h}$.
Then $h$ is extended to a linear automorphism of $V_{L^\circ}$ mapping $V_{\alpha+L}$ to $V_{\bar{g}(\alpha)+L}$ as in (\ref{act}).
In particular, {its restriction to $V_L$ is an element of $\Aut(V_L)$.}
By Theorem \ref{PAVL}, there exists an element $f\in N(V_L)$ such that $g=f\circ h_{|V_L}\in \AutO(\hat{L})$.
By Proposition \ref{RAM} (3), we have $V_{\alpha+L}^f\cong V_{\alpha+L}$.
Hence we obtain $V_{\alpha+L}^g\cong h(V_{\alpha+L})=V_{\bar{g}(\alpha)+L}$.
In particular, we have $[\mu]^g=[\bar{g}(\mu)]$ for $\mu\in L^\circ\setminus (L/2)$, $\{[\lambda]^{\pm,g}\}=\{[\bar{g}(\lambda)]^\pm\}$ for $\lambda\in L^\circ\cap (L/2)$ and $[0]^{\pm,g}=[0]^\pm$.

We assume that $\bar{g}=1$.
Then $g=f_{\beta}$ for some $\beta\in L^\circ$.
Let $\lambda$ be an element of $L^\circ\cap(L/2)$.
Let $\psi_{\beta}:V_{\lambda+L}\to V_{\lambda+L}$ be the isomorphism of $V_L$-modules defined by $x\ot e(\lambda+\gamma)\mapsto (\sqrt{-1})^{\langle \beta,2\lambda\rangle}x\ot e(\lambda+\gamma)$, $\gamma\in L$.
Then $Y(gv,z)\psi_\beta w=\psi_\beta Y(v,z)w$ for $v\in V_L$, $w\in V_{\lambda+L}$.
It is easy to see that $\theta_{V_L}\psi_\beta=(-1)^{\langle \beta,2\lambda\rangle}\psi_\beta\theta_{V_L}$.
Therefore we obtain this proposition.
\qe

The following lemma is useful in calculating the automorphism group of $V_L^+$.

\bl{PStab} Let $N$ be a sublattice of $L^\circ\cap (L/2)$ such that $L\subset N$.
Set 
\begin{eqnarray*}
G=\{g\in O(\hat{L})|\ [\lambda]^{\varepsilon,g}=[\lambda]^\varepsilon\ {\rm for}\ {}^\forall\varepsilon\in\{\pm\},\ {}^\forall\lambda\in N\}.
\end{eqnarray*}
Then the following sequence
\begin{eqnarray*}
1\to\{f\in\Hom(L,\Z_2)|\ f(2N)=0\}\hookrightarrow G\ \bar{\to}\ \{g\in O(L)|\ g=1\ {\rm on}\ 2N/2L\}\to1
\end{eqnarray*}
is exact.
\el
\proof\ By the exact sequence (\ref{exactauto}), it suffices to determine both $\Hom(L,\Z_2)\cap G$ and $\overline{G}$.
By Proposition \ref{actioncoset}, we have $\Hom(L,\Z_2)\cap G=\{f\in\Hom(L,\Z_2)|\ f(2W)=0\}$.
Let us determine $\overline{G}$.
By Proposition \ref{actioncoset}, $\overline{G}$ acts trivially on $2N/2L$.
Hence we have $\overline{G}\subseteq G_0=\{g\in O(L)|\ g=1\ {\rm on}\ 2N/2L\}$.
Conversely, let $h$ be an element of $G_0$.
Then there exists $g\in O(\hat{L})$ such that $\bar{g}=h$ by the exact sequence (\ref{exactauto}).
Let $s:2N/2L\to\Z_2$ be the map defined by $g(a)=(-1)^{s(\bar{a})}a$ for $a\in\hat{2N}$.
Then $s$ is a homomorphism {of groups}.
Since $2N/2L$ is a subspace of $L/2L$ over $\F_2$, $s$ is extended to $s_0\in\Hom(L,\Z_2)$ such that $s_0=s$ on $2N/2L$.
Hence we have an element $s_0g\in G$ such that $\overline{s_0g}=h$.
Therefore we obtain $\overline{{G}}=G_0$.
\qe

\subsection{Extra automorphisms of $V_L^+$}\label{SEAVL+}
In this subsection, we construct some automorphism of $V_L^+$ called an {\it extra automorphism} for an even lattice $L$ obtained by Construction B along the line of \cite{FLM2}.
{The term ``extra automorphism"} means that it does not fix the isomorphism class of $V_L^-$.
We will show in Proposition \ref{CMain} that an automorphism of $V_L^+$ is extra if and only if it is not obtained as the restriction of an automorphism of $V_L$ to $V_L^+$.

Let $C$ be a doubly even code of length $n$ and let $\{\alpha_i|\ i\in\Omega_n\}$ be an orthogonal basis of $\R^n$ of norm $2$.
Let $L_B(C)$ be the lattice obtained by Construction B from $C$ associated with $\{\alpha_i|\ i\in\Omega_n\}$.
Then $L_B(C)$ is an even lattice of rank $n$ by Proposition \ref{PCOB} (1).
Set $L_A(C)=L_B(C)+\Z\alpha_1$, and fix $a_k\in\hat{L}_A(C)$ for each $k\in\Omega_{n}$ such that $\bar{a}_k=\alpha_k$.
Let $\sigma$ be the operator of $V_{L_A(C)}$ defined by
\begin{eqnarray}
\sigma&=&\prod_{k=1}^n\exp{((1+\sqrt{-2})(a_k)_0)}\exp(\sqrt{\frac{-1}{2}}(a_k^{-1})_0)\exp{((-1+\sqrt{-2})(a_k)_0)}.\label{extra}
\end{eqnarray}
Then $\sigma$ is an automorphism of the VOA $V_{L_A(C)}$.

\br{RSI}{\rm The automorphism $\sigma$ depends on the choice of the orthogonal basis $\{\alpha_i|\ i\in\Omega_n\}$ of norm $2$ which is used in the construction of the lattice $L_B(C)$.}\er

The following proposition is a slight modification of \cite[Proposition 12.2.5]{FLM2}.
\bp{PGVL+} Let $J$ be a set of generators of an even lattice $L$ of rank $n$ and let $\{h_k|\ k\in\Omega_n\}$ be a basis of $\h_L$.
Then the VOA $V_L^+$ is generated by $\{h_i(-1)h_j(-1),\ a^+|\ i,j\in\Omega_n\ a\in\hat{J}\}$, where $a^+=a+\theta_{V_L}(a)$.
\ep

Let $J$ be a set of generators of $L_B(C)$ given in Proposition \ref{PGCB}.
By the proposition above, $\tilde{J}=\{\alpha_i(-1)\alpha_j(-1),\ a^+|\ i,j\in\Omega_n\ a\in\hat{J}\}$ is a set of generators of $V_{L_B(C)}^+$.
The following proposition follows from \cite[Theorem 11.2.1, Corollary 11.2.4 and Proposition 12.2.1]{FLM2}.
\bp{PAS}{\rm \cite{FLM2}} The automorphism $\sigma$ preserves the set $\tilde{J}$.
In particular $\sigma$ is an automorphism of $V_{L_B(C)}^+$.
Moreover the order of $\sigma$ is two.
\ep

The following proposition is a slight generalization of \cite[Proposition 12.2.7]{FLM2}.
\bp{PPSVL+} The automorphism $\sigma$ sends $V_{L_B(C)}^-$ and $V_{\alpha+L_B(C)}^-$ to $V_{\alpha_1+L_B(C)}^+$ and $V_{\alpha_1+L_B(C)}^-$ respectively.
In particular, $[0]^{-,\sigma}= [\alpha_1]^+$ and $[\alpha_1]^{-,\sigma}=[\alpha_1]^-$.
\ep
This proposition shows that $\sigma$ is an extra automorphism of $V_{L_B(C)}^+$.

\br{RFLMEXT}{\rm In \cite{FLM2}, it is assumed that the all-one codeword belongs to $C$.
However, this assumption is not necessary to obtain the extra automorphism $\sigma$ of $V_{L_B(C)}^+$.}
\er

\section{Automorphism group of $V_L^+$}\label{ChAGVL+}
In this section, we study the automorphism group $\Aut(V_L^+)$ of $V_L^+$ by using its action on isomorphism classes of irreducible $V_L^+$-modules, and give a method of determining $\Aut(V_L^+)$.
\subsection{VOA fixed by a finite abelian group}\label{SAVGA}
In this subsection, we study the automorphism groups of a simple VOA and its subVOA obtained as the fixed-point subspace under an abelian automorphism group.
We will apply the results of this section to the VOAs $V_L$ and $V_L^+$ in Section \ref{SAGVL+}.

Let $V$ be a simple VOA and let $D$ be a finite abelian subgroup of $\Aut(V)$.
For $d\in \Hom(D,\C^\times)$ we set $V(d)=\{v\in V|\ b(v)=d(b)v\ \forall b\in D\}$.
Since the field $\C$ is algebraically closed and $D$ is finite abelian, we have $V=\op_{d\in \Hom(D,\C^\times)}V(d)$.
Let $1_D$ denote the identity element of $\Hom(D,\C^\times)$.
Then $V(1_D)$ is a subVOA of $V$ and $V(d)$ is a $V(1_D)$-module.
Clearly the restriction $Y_{d_1,d_2}(u,z)$ of $Y(u,z)$ to $V(d_2)$ for $u\in V(d_1)$ is an intertwining operator of type $\left(\IntType{V(d_1)}{V(d_2)}{V(d_1d_2)}\right)$.
In this paper, we will use the following properties of $V(d)$ studied in \cite{DM1}.

\bt{PSV}{\rm \cite{DM1}} Let $V$ be a simple VOA and let $D$ be a finite abelian subgroup of $\Aut(V)$.
Then $\{V(d)|\ d\in \Hom(D,\C^\times)\}$ is a set of non-isomorphic irreducible $V(1_D)$-modules.
In particular, $V(1_D)$ is a simple VOA.
\et

The following proposition obtained in \cite{DM2} is crucial for our purpose.
\bp{PDM2}{\rm \cite[Proposition 5.3]{DM2}} Let $V=(V,Y,\1,\w)$ be a simple VOA and let $D$ be a finite abelian subgroup of $\Aut(V)$.
Suppose that $V(d_1)\times V(d_2)=V(d_1d_2)$ for all $d_i\in \Hom(D,\C^\times)$.
If $\tilde{V}=(V,\tilde{Y},\1,\w)$ is also a simple VOA with $Y(v,z)=\tilde{Y}(v,z)$ for all $v\in V(1_D)$ then there exists an isomorphism $\tau:\tilde{V}\to V$ of VOAs such that $\tau_{|V(1_D)}$ is the identity operator and $\tau$ preserves $V(d)$ for all $d\in \Hom(D,\C^\times)$.
\ep

Clearly, the normalizer $N_{\Aut(V)}(D)$ of $D$ in $\Aut(V)$ preserves $V(1_D)$ and $\{V(d)|\ d\in \Hom(D,\C^\times)\}$ as a set of $V(1_D)$-modules.
We note that if $g\in N_{\Aut(V)}(D)$ satisfies $g(V(d_1))=V(d_2)$ then $V(d_1)^{g_{|V(1_D)}}\cong V(d_2)$ as $V(1_D)$-modules.
Consider the restriction homomorphism
\begin{eqnarray*}
\varphi_V:N_{\Aut(V)}(D)\to H,
\end{eqnarray*}
where 
\begin{eqnarray}
H=\Bigl\{h\in\Aut(V(1_D))\Big\arrowvert\ \{V(d)^h|\ d\in \Hom(D,\C^\times)\}=\{V(d)|\ d\in \Hom(D,\C^\times)\}\Bigr\}.\label{Def:H}
\end{eqnarray}

\bt{sc} Let $V=(V,Y,\1,\w)$ be a simple VOA and let $D$ be a finite abelian subgroup of $\Aut(V)$.
Suppose that $V(d_1)\times V(d_2)=V(d_1d_2)$ for all $d_i\in \Hom(D,\C^\times)$.
Then the restriction homomorphism $\varphi_V$ is surjective and $\Ker\varphi_V=D$.
\et
\proof\ First we show the surjectivity of $\varphi_V$.
Let $g$ be an element of $H$.
For any $d_1\in \Hom(D,\C^\times)$, there exists a unique element $d_2$ of $\Hom(D,\C^\times)$ such that $V(d_2)\cong V(d_1)^g$ as $V(1_D)$-modules by Proposition \ref{PSV} (2), and let $\psi_{d_2}:V(d_2)\to V(d_1)^g$ be an isomorphism of $V(1_D)$-modules.
Then
\begin{eqnarray}
Y_{1_D,d_1}(gv,z)\psi_{d_2}=\psi_{d_2} Y_{1_D,d_2}(v,z)\label{Eq:Const1}
\end{eqnarray}
for $v\in V(1_D)$.
In particular, we may choose 
\begin{eqnarray}
\psi_{1_D}=g\label{Eq:Const2}
\end{eqnarray}
since $V(1_D)^g\cong V(1_D)$ and $Y_{{1_D},{1_D}}(gv,z)g=gY_{{1_D},{1_D}}(v,z)$ for $v\in V({1_D})$.
Then we have a linear automorphism $\psi=\op_{d\in \Hom(D,\C^\times)}\psi_d$ of $V$.
Set $Y^g(v,z)=\psi^{-1} Y(\psi(v),z)\psi$ for $v\in V$.
By (\ref{Eq:Const2}), it is easy to see that $\tilde{V}=(V,Y^g,\1,\w)$ has a VOA structure.
Moreover, for $v\in V$
\begin{eqnarray*}
\psi Y^g(v,z)=\psi\psi^{-1} Y(\psi(v),z)\psi=Y(\psi (v),z)\psi.\label{Eq:Const3}
\end{eqnarray*}
Hence $\psi$ is an isomorphism of VOAs from $V$ to $\tilde{V}$.
By (\ref{Eq:Const1}) and (\ref{Eq:Const2}), $Y^g(v,z)=Y(v,z)$ for $v\in V({1_D})$.
{Applying Proposition \ref{PDM2} to $V$ and $\tilde{V}$, we obtain an isomorphism $\tau:\tilde{V}\to V$ such that $\tau_{|V({1_D})}$ is the identity operator and {$\tau$} preserves $V(d)$ for all $d\in \Hom(D,\C^\times)$.
Then $\tau\circ\psi$ is an automorphism of the VOA $V$ such that its restriction to $V({1_D})$ is $g$ and $\tau\circ\psi\in N_{\Aut(V)}(D)$.
Therefore $\varphi_V(\tau\circ\psi)=g$, and $\varphi_V$ is surjective.

Next we determine the kernel of $\varphi_V$.
Let $h$ be an element of $\Ker\varphi_V$.
By Schur's lemma, $h$ acts on $V(d)$ by a scalar $\lambda_d\in\C^\times$ for each $d\in \Hom(D,\C^\times)$.
In particular, $\lambda_{1_D}=1$.
By the assumptions on the fusion rules, we have $\lambda_{d_1}\lambda_{d_2}=\lambda_{d_1d_2}$ for $d_i\in \Hom(D,\C^\times)$.
Hence $h\in D$.
Clearly $D\subseteq\Ker\varphi_V$.
Therefore we obtain $\Ker\varphi_V=D$.
\qe

\subsection{$V_L^+$ for a $2$-elementary totally even lattice}\label{SVL+2ETE}
In this subsection, we study the irreducible $V_L^+$-modules and their fusion rules for a $2$-elementary totally even lattice $L$.
The results of this section will be used for the description of $\Aut(V_L^+)$ when $L$ is isomorphic to $\sqrt2E_8$ and the Barnes-Wall lattice of rank $16$ in Section \ref{SsSqrt2R} and \ref{SsBW}.

Let $L$ be a $2$-elementary totally even lattice $L$ of rank $n$, namely $2L^\circ\subset L$ and both $L$ and $\sqrt2L^\circ$ are even.
Then $L^\circ/L\cong\Z_2^m$ for some $m\in\Z$.
We note that the determinant of $L$ is equal to $2^m$.
By Theorem \ref{Eq:CVL+} and Remark \ref{Risomod}, $V_L^+$ has exactly $2^{m+1}$ non-isomorphic irreducible modules $V_{\lambda+L}^\pm$ ($\lambda+L\in L^\circ/L$) of untwisted type.
Let us consider irreducible $V_L^+$-modules of twisted type.
By Proposition \ref{PTC}, such a $V_L^+$-module is uniquely determined by an element of $X(Z(\hat{L}/K_L))$.
Since $2L^\circ\subset L$, we have $Z(\hat{L}/K_L)=\hat{2L^\circ}/K_L$.
Since $\langle \alpha,\alpha\rangle\in4\Z$ for all $\alpha\in 2L^\circ$, the group $\hat{2L^\circ}/K_L$ is an elementary abelian $2$-group.
We choose a complement $G$ to $\langle\kappa_LK_L\rangle$ in $\hat{2L^\circ}/K_L$, which is a subgroup of $\hat{2L^\circ}/K_L$ such that $\hat{2L^\circ}/K_L=G\langle\kappa_LK_L\rangle$ and $G\cap\langle\kappa_LK_L\rangle=\{K_L\}$.
Then the map $G\to 2L^\circ/2L$, $a\mapsto \bar{a}$ is an isomorphism of groups, so we identify $G$ with $2L^\circ/2L$.
Let $\chi\in X(\hat{2L^\circ}/K_L)$.
Then there exists a unique coset $\lambda+L\in L^\circ/L$ such that $\chi(\beta)=(-1)^{\langle \lambda,\beta\rangle}$ for $\beta\in 2L^\circ/2L$.
We denote such a character by $\chi_\lambda$.
We note that this notation depends on the choice of the complement.
Hence $V_L^+$ has exactly $2^{m+1}$ non-isomorphic irreducible modules $V_L^{T_{\chi_\lambda},\pm}$ ($\lambda+L\in L^\circ/L$) of twisted type.
Therefore $V_L^+$ has exactly $2^{m+2}$ non-isomorphic irreducible modules when $L$ is a $2$-elementary totally even lattice with determinant $2^m$.

Let $S_L$ be the set of all isomorphism classes of irreducible $V_L^+$-modules.
By the arguments above, $S_L=\{[\lambda]^\pm,[\chi_\lambda]^\pm|\ \lambda+L\in L^\circ/L\}$.
Let us describe the fusion rules of $V_L^+$.
The fusion rules in this case is clearer than those in general case.
We regard $\{\pm\}$ as the group of order $2$ generated by $-$.
Let $\nu:L^\circ\to\{\pm\}$ be the map defined by $\nu(\alpha)=+$ if $\langle \alpha,\alpha\rangle\in2\Z$, and $\nu(\alpha)=-$ if $\langle \alpha,\alpha\rangle\in 1+2\Z$.
Applying the results of \cite{ADL} to our case, we obtain the following proposition.

\bp{fusionoflattice}
{Let $L$ be a $2$-elementary totally even lattice.
Then the fusion rules of $V_L^+$ are described as follows:}
\eqn
[\lambda_1]^\delta\times[\lambda_2]^\varepsilon&=&[\lambda_1+\lambda_2]^{\delta\varepsilon},\\
{} [\lambda_1]^\delta\times[\chi_{\lambda_2}]^\varepsilon&=&[\chi_{\lambda_1+\lambda_2}]^{\delta\varepsilon\nu(\lambda_2)\nu(\lambda_1+\lambda_2)},\\
{}[\chi_{\lambda_1} ]^\delta\times[\chi_{\lambda_2}]^\varepsilon&=&[\lambda_1+\lambda_2]^{\delta\varepsilon\nu(\lambda_1)\nu(\lambda_2)},
\eeqn
where $\delta,\varepsilon\in\{\pm\}\cong\Z_2$ and $\lambda_1,\lambda_2\in L^\circ$.
\ep

Now, we view the formal product $\times$ as a binary operation on $S_L$.
By Proposition \ref{PFHL} and \ref{RFusionFo} (3), this operation is associative and commutative.
It is easy to see that $W\times W=[0]^+$ for all $W\in S_L$.
Hence $S_L$ has a structure of an elementary abelian $2$-group.
So the set $S_L$ forms an $(m+2)$-dimensional vector space over $\F_2$ under the fusion rules.

\br{ordoffusion} {\rm Let $L$ be an even lattice.
{If $2L^\circ \not\subset L$ then $[\mu]\times[\mu]=[0]+[2\mu]$ for $\mu\in L^\circ\setminus (L/2)$.
If $\sqrt2(L^\circ\cap(L/2))$ is not even then $[\lambda]^+\times [\lambda]^+=[0]^-$ for $\lambda\in L^\circ\cap (L/2)$ with $\langle\lambda,\lambda\rangle\in 1/2+\Z$.
Hence the set $S_L$ forms a vector space over $\F_2$ under the fusion rules if and only if $L$ is $2$-elementary totally even.}}
\er

Let us consider the action of $\AutO(\hat{L})$ on $S_L$.

\bl{actioncoset2} For $f_\beta\in\Hom(L,\Z_2)$, $\lambda\in L^\circ$ and $\varepsilon\in\{\pm\}$, we have
\eqn
[\chi_\lambda]^{\varepsilon,f_\beta}=[\chi_{\lambda+\beta}]^\varepsilon.
\eeqn
In particular, $f_\beta$ fixes all isomorphism classes of irreducible $V_L^+$-modules of twisted type if and only if $\beta\in L/2L^\circ$.
\el
\proof\ We set $\chi_\lambda^{(\beta)}(a)=(-1)^{\langle\beta,\bar{a}\rangle}\chi_\lambda(a)$, $a\in \hat{L}$.
By the definition of $\chi_\lambda$, we have $\chi_\lambda^{(\beta)}=\chi_{\lambda+\beta}$.
It is easy to see that the actions of $\Cent(\hat{L}/K_L)$ on $[\chi_\lambda]^{\varepsilon,f_\beta}$ and $[\chi_\lambda^{(\beta)}]^\varepsilon$ are the same, hence $[\chi_\lambda]^{\varepsilon,f_\beta}=[\chi^{(\beta)}]^\varepsilon$.
Thus we obtain $[\chi_\lambda]^{\varepsilon,f_\beta}=[\chi_{\lambda+\beta}]^\varepsilon$.\qe

\bl{kernelofauto3} Let $F$ be the group consisting of all elements of $\AutO(\hat{L})/\langle\theta_{V_L}\rangle$ that fix all isomorphism classes of irreducible $V_L^+$-modules.
Then the sequence
\begin{eqnarray*}
1\to\{f_\beta|\ \beta\in L/2L^\circ\}\hookrightarrow F\ \bar{\to}\ \{g\in O(L)|\ g=1\ {\rm on}\ 2L^\circ/2L\}\to1
\end{eqnarray*}
is exact.
\el
\proof\ By Lemma \ref{PStab} and Lemma \ref{actioncoset2}, we have $F\cap\Hom(L,\Z_2)=\{f_\beta|\ \beta\in L/2L^\circ\}$ and $\bar{F}\subseteq\{g\in O(L)|\ g=1\ {\rm on}\ 2L^\circ/2L\}$.
Let us show that this inclusion becomes equality.
Let $h_0$ be an automorphism of $L$ acting trivially on $2L^\circ/2L$ and let $h$ be an element of $\AutO(\hat{L})$ such that $\bar{h}=h_0$.
Then there exists a map $s:2L^\circ/2L\to\Z_2$ such that $h(a)=(-1)^{s(\bar{a})}a$ for $a\in \hat{2L^\circ}$.
It is easy to see that $s$ is a {group} homomorphism and $s=f_\beta$ for some $\beta\in L^\circ/2L^\circ$.
Hence $f_\beta h$ fixes all isomorphism classes of irreducible modules of twisted type.
By Theorem \ref{fusionoflattice}, $fh_\beta$ also fixes all isomorphism class of irreducible $V_L^+$-modules of untwisted type.
Hence $fh_\beta\in F$ and $\overline{f_\beta h}=h_0$.
Thus we obtain $\bar{F}=\{g\in O(L)|\ g=1\ {\rm on}\ 2L^\circ/2L\}$.\qe

We saw that the set $S_L$ of all isomorphism classes of irreducible $V_L^+$-modules forms a vector space over $\F_2$ under the fusion rules.
We will introduce a quadratic form on $S_L$ preserved by the action of $\Aut(V_L^+)$.
We are interested in the cases $n\in8\Z$ since we will consider $2$-elementary totally even lattices $\sqrt2E_8$ and $\BW$ later.
So, in the rest of this subsection, we assume $n\in8\Z$.
Then by Proposition \ref{PCH}, $\ch W$ belongs to either $\Z[q]$ or $q^{1/2}\Z[q]$ for each $W\in S_L$.

{Let $q_L:S_L\to\F_2$ be the map defined by}
\begin{eqnarray*}
q_L(W)=\left\{\begin{array}{cl}
 \mbox{$0$} & \mbox{${\rm if}\ \ch(W)\in\Z[q]$},\\
 \mbox{$1$} & \mbox{${\rm if}\ \ch(W)\in q^{1/2}\Z[q]$},
\end{array}
\right.
\end{eqnarray*}
where $W\in S_L$.
In the other word,
\eqn
q_L([\lambda]^\pm)&=&\langle\lambda,\lambda\rangle\mod2,\\
q_L([\chi_\lambda]^+)&=&\left\{\begin{array}{cl}
 \mbox{$1$,} & \mbox{$(n\in8+16\Z)$},\\
 \mbox{$0$,} & \mbox{$(n\in16\Z)$},
\end{array}
\right.\\
q_L([\chi_\lambda]^-)&=&\left\{\begin{array}{cl}
 \mbox{$0$,} & \mbox{$(n\in8+16\Z)$},\\
 \mbox{$1$,} & \mbox{$(n\in16\Z)$},
\end{array}
\right.
\eeqn
where $\lambda\in L^\circ$.
By Proposition \ref{RAM} (2), $q_L(W)=q_L(W^g)$ for $W\in S_L$, $g\in\Aut(V_L^+)$.

\bt{PQuad} Let $L$ be a $2$-elementary totally even lattice whose rank is a multiple of $8$.
Then the map $q_L$ is a quadratic form associated with a non-singular symplectic form on $S_L$ such that $q_L(W)=q_L(W^g)$ for all $W\in S_L$, $g\in\Aut(V_L^+)$.
\et
\proof\ Let us check that the form $(\cdot,\cdot)_L:S_L\times S_L\to\F_2$ defined by $(x,y)_L=q_L(x)+q_L(y)-q_L(x+y)$ is non-singular symplectic.
Let $\pi:\{\pm\}\to\F_2$ be an isomorphism of groups such that $\pi(+)=0$ and $\pi(-)=1$.
Let $\lambda_1,\lambda_2$ be vectors in $L^\circ$ and let $\delta,\varepsilon$ be elements of $\{\pm\}$.
It is easy to see that
\eqn
q_L([\lambda_1]^\delta)+q_L([\lambda_2]^\varepsilon)+q_L([\lambda_1+\lambda_2]^{\delta\varepsilon})&=&\langle\lambda_1,\lambda_1\K+\langle\lambda_2,\lambda_2\K+\langle\lambda_1+\lambda_2,\lambda_1+\lambda_2\K\\{} &=&2\langle\lambda_1,\lambda_2\K,\\{}
q_L([\lambda_1]^\delta)+q_L([\chi_{\lambda_2}]^\varepsilon)+q_L([\chi_{\lambda_1+\lambda_2}]^{\delta\varepsilon\nu(\lambda_2)\nu(\lambda_1+\lambda_2)})&=&\langle\lambda_1,\lambda_1\K+\pi(\varepsilon)+\pi(\delta\varepsilon)\\{}&&+\langle\lambda_2,\lambda_2\K+\langle\lambda_1+\lambda_2,\lambda_1+\lambda_2\rangle\\{}
&=&2\langle\lambda_1,\lambda_2\K+\pi(\delta),\\{}
q_L([\chi_{\lambda_1}]^\delta)+q_L([\chi_{\lambda_2}]^\varepsilon)+q_L([\lambda_1+\lambda_2]^{\delta\varepsilon\nu(\lambda_1)\nu(\lambda_2)})&=&\pi(\delta)+\pi(\varepsilon)+\langle\lambda_1+\lambda_2,\lambda_1+\lambda_2\rangle\\{}
&=&\pi(\delta\varepsilon)+\langle\lambda_1+\lambda_2,\lambda_1+\lambda_2\K.
\eeqn
Hence we have
\eqn
([{\lambda_1}]^\delta,[{\lambda_2}]^\varepsilon)_L&=&2\langle{\lambda_1},{\lambda_2}\K,\\{}
([{\lambda_1}]^\delta,[\chi_{\lambda_2}]^\varepsilon)_L&=&2\langle{\lambda_1},{\lambda_2}\K+\pi(\delta),\\{}
([\chi_{\lambda_1}]^\delta,[\chi_{\lambda_2}]^\varepsilon)_L&=&\langle{\lambda_1}+{\lambda_2},{\lambda_1}+{\lambda_2}\K+\pi(\delta\varepsilon).
\eeqn
We will show that $(\cdot,\cdot)_L$ coincides with the non-singular symplectic form $(\cdot,\cdot)$ defined by
\eqn
([{\lambda_1}]^+,[{\lambda_2}]^+)=2\B{\lambda_1},{\lambda_2}\K,\ ([{\lambda_1}]^+,[\chi_0]^\pm)=0,\ ([\chi_0]^+,[\chi_0]^-)=1,\ ([\chi_0]^\pm,[\chi_0]^\pm)=0.
\eeqn
It is easy to see that
\eqn
([{\lambda_1}]^\delta,[{\lambda_2}]^\varepsilon)&=&([{\lambda_1}]^+,[{\lambda_2}]^+)=2\langle{\lambda_1},{\lambda_2}\K,\\{}
([{\lambda_1}]^\delta,[\chi_{\lambda_2}]^\varepsilon)&=&([{\lambda_1}]^\delta,[{\lambda_2}]^{\varepsilon\nu({\lambda_2})})+([{\lambda_1}]^\delta,[\chi_0]^+)=2\langle{\lambda_1},{\lambda_2}\K+\pi(\delta),\\{}
([\chi_{\lambda_1}]^\delta,[\chi_{\lambda_2}]^\varepsilon)&=&([{\lambda_1}]^{\delta\nu({\lambda_1})},[{\lambda_2}]^{\varepsilon\nu({\lambda_2})})+([{\lambda_1}]^{\delta\nu({\lambda_1})},[\chi_0]^+)\\{}
& &+([\chi_0]^+,[{\lambda_2}]^{\varepsilon\nu({\lambda_2})})+([\chi_0]^+,[\chi_0]^+)\\{}
&=&2\langle{\lambda_1},{\lambda_2}\K+\pi(\delta\nu({\lambda_1}))+\pi(\varepsilon\nu({\lambda_2}))\\{}
&=&\langle{\lambda_1}+{\lambda_2},{\lambda_1}+{\lambda_2}\K+\pi(\delta\varepsilon).
\eeqn
Hence we have $(\cdot,\cdot)_L=(\cdot,\cdot)$.
In particular, the form $(\cdot,\cdot)_L$ is non-singular symplectic.\qe

We consider the type of the quadratic from $q_L$.
By the definition of $q_L$, $S_L$ decomposes into an orthogonal direct sum $S_L=S_0\op\{[0]^\pm,[\chi_0]^\pm\}$ with respect to the symplectic form associated with $q_L$, where $S_0=\{[\lambda]^+|\ \lambda\in L^\circ/L\}$ is a subspace of $S_L$.
On the other hand, $L^\circ/L$ is a vector space over $\F_2$ equipped with a non-singular quadratic form $q(\lambda+L)=\langle\lambda,\lambda\rangle\mod2$, $\lambda\in L^\circ/L$.
The linear isomorphism $f:L^\circ/L\to S_0$, $\lambda+L\mapsto [\lambda]^+$ show that the types of the restriction of $q_{L}$ to $S_0$ and $q$ are the same.
Since the type of the restriction of $q_L$ to $\{[0]^\pm,[\chi_0]^\pm\}$ is plus, the type of $q_L$ is equal to that of $q$.
Hence we obtain the following remark.
\br{symonlattice}{\rm The type of $q_L$ is equal to that of the quadratic form $q$ on $L^\circ/L$ given by $q(\lambda+L)=\langle\lambda,\lambda\rangle\mod2$, $\lambda\in L^\circ/L$.}\er

\subsection{Action of the automorphism group of $V_L^+$ on isomorphism classes of the irreducible modules}\label{SAGVL+}
In this subsection, we study the action of $\Aut(V_L^+)$ on the set of isomorphism classes of irreducible $V_L^+$-modules.
Among other things, we consider the stabilizer of the isomorphism class $[0]^-$ of $V_L^-$ and the orbit of $[0]^-$.

Let $L$ be an even lattice of rank $n$ and let $S_L$ be the set of all isomorphism classes of irreducible $V_L^+$-modules.
Let $H_L$ be the stabilizer of $[0]^-$ under the action of $\Aut(V_L^+)$ on $S_L$, that is, $H_L=\{g\in\Aut(V_L^+)|\ [0]^{-,g}= [0]^-\}$.

First we determine the group $H_L$.
It is easy to see that $V_L$ and $\langle\theta_{V_L}\rangle$ satisfy the assumption of Theorem \ref{sc}.
We note that $C_{\Aut(V_L)}(\theta_{V_L})=N_{\Aut(V_L)}(\langle\theta_{V_L}\rangle)$ since the order of $\theta_{V_L}$ is two.
By Theorem \ref{sc}, the restriction homomorphism from $C_{\Aut(V_L)}(\theta_{V_L})$ to $H_L$ is surjective and its kernel is $\langle\theta_{V_L}\rangle$.}
Hence we have the following proposition.

\bp{CMain} The stabilizer $H_L$ of $[0]^-$ is $C_{\Aut(V_L)}(\theta_{V_L})/\langle\theta_{V_L}\rangle$.
\ep

As an application of Proposition \ref{CMain}, {we have the following proposition.}

\bp{co4} The group $\Aut(V_L^+)$ is finite if and only if $L$ has no roots.
\ep
\proof\ Suppose that $L$ has roots, and let $\alpha$ be a root of $L$.
Then $v=e(\alpha)+\theta_{V_L}(e(\alpha))$ is a vector in $(V_L^+)_1$.
It is easy to see that $v_0(e(2\alpha)+\theta_{V_L}(e(2\alpha)))\neq 0$, which implies that $\exp(v_0)$ generates an infinite subgroup of $\Aut(V_L^+)$.

Conversely, suppose that $L$ has no roots.
Then $C_{\Aut(V_L)}(\theta_{V_L})$ is finite by Proposition \ref{co1}.
By Theorem \ref{Eq:CVL+} and Remark \ref{Risomod}, the cardinality of $S_L$ is finite.
Hence $|\Aut(V_L^+)|\le |S_L|\cdot|C_{\Aut(V_L)}(\theta_{V_L})/\langle\theta_{V_L}\rangle|<\infty$.\qe

Set $\Orbit=\{W\in S_L|\ W\times W=[0]^+,\ \ch W=\ch [0]^-\}$.
Then we have the following lemma.

\bl{orbit} Let $L$ be an even lattice of rank $n$.
\begin{enumerate}
\item If $n\neq8,16$ then $\Orbit\subseteq\{[0]^-,[\lambda]^\pm|\ \lambda\in L^\circ\cap (L/2),\ |(\lambda+L)_2|=2n+|L_2|\}$.
\item If $n=8$ then $\Orbit\subseteq\{[0]^-,[\lambda]^\pm,\ [\chi]^-|\ \lambda\in L^\circ\cap (L/2),\ |(\lambda+L)_2|=2n+|L_2|,\ \chi\in X(\Cent(\hat{L}/K_L))\}$.
\item If $n=16$ then $\Orbit\subseteq\{[0]^-,[\lambda]^\pm,\ [\chi]^+|\ \lambda\in L^\circ\cap (L/2),\ |(\lambda+L)_2|=2n+|L_2|,\ \chi\in X(\Cent(\hat{L}/K_L))\}$.
 \end{enumerate}
\el
\proof\ Proposition \ref{RFusionFo} (1) shows that $[\mu]\notin \Orbit$ for $\mu\in L^\circ\setminus(L/2)$.
Let $\lambda$ be a vector in $L^\circ\cap (L/2)$ and let $\varepsilon$ be an element of $\{\pm\}$.
If $[\lambda]^\varepsilon\in\Orbit$ then $\dim (V_{\lambda+L}^\varepsilon)_1=\dim (V_L^-)_1$.
By Proposition \ref{PCH}, we have $|(\lambda+L)_2|=2n+|L_2|$.
Let $\chi$ be an element of $X(\Cent(\hat{L}/K_L))$.
By Proposition \ref{PCH}, if $[\chi]^+\in \Orbit$ then $n=16$, and if $[\chi]^-\in \Orbit$ then $n=8$.
\qe

Next we consider the orbit $Q_L$ of the isomorphism class $[0]^-$ of $V_L^-$ or $Q_L=\{W\in S_L|\ W^g=[0]^-,\ {}^\exists g\in\Aut(V_L^+)\}$.
Let us study an even lattice $L$ satisfying the following condition:
\begin{enumerate}[(I)]
\item $Q_L$ contains isomorphism classes of irreducible $V_L^+$-modules of twisted type.\label{Assumption1}
\end{enumerate}

\bl{orbit2} Let $L$ be an even lattice satisfying (\ref{Assumption1}).
Then $2L^\circ\subset L$.
\el
\proof\ Suppose that $2L^\circ\not\subset L$ and let $\mu\in L^\circ\setminus (L/2)$.
By Proposition \ref{RFusionFo} (1), the fusion rule $[0]^-\times [\mu]=[\mu]$ holds.
By the condition (\ref{Assumption1}), there exists $g\in\Aut(V_L^+)$ such that $[0]^{-,g}$ is of twisted type.
Then the fusion rule $[0]^{-,g}\times[\mu]^g=[\mu]^g$ contradicts Proposition \ref{RFusionFo} (2).
\qe

The identity relating the theta series of a lattice and the Dedekind $\eta$-function is crucial in the following proposition.

\bp{orbit3} Let $L$ be an even lattice satisfying (\ref{Assumption1}).
Then $L$ is a $2$-elementary totally even lattice of rank $8$ or $16$.
Moreover, if $L$ does not have roots then $L$ is isomorphic to either $\sqrt2E_8$ or the Barnes-Wall lattice $\BW$ of rank $16$.
\ep
\proof\ By Lemma \ref{orbit}, the rank $n$ of $L$ is $8$ or $16$.
By Lemma \ref{orbit2}, $2L^\circ\subset L$.
Let $m$ be the integer given by $|L^\circ/L|= 2^m$.
We note that $2^m$ is the determinant of $L$.

Let us check that $\langle v,v\rangle\in\Z$ for all $v\in L^\circ$.
Set $\eps=+$ if $n=16$, and $\eps=-$ if $n=8$.
Let $\chi$ be an element of $X(Z(\hat{L}/K_L))$ such that $[0]^{-,g}=[\chi]^{\eps}$ for some $g\in\Aut(V_L^+)$.
By Theorem \ref{PTC}, the dimension of $T_{\chi}$ is equal to $|L/(L\cap 2L^\circ)|^{1/2}=2^{k/2}$, where $k=n-m$.
By the equation $\ch [0]^-=\ch [\chi]^\eps$ and Proposition \ref{PCH}, the theta series of $L$ is described as follows:
\begin{eqnarray}
n=8:\ \Theta_L(q)&=&\frac{\eta(\tau)^{16}}{\eta(2\tau)^8}+2^{\frac{k}{2}}\frac{\eta(\tau)^{16}}{\eta(\tau/2)^8} -2^{\frac{k}{2}}\frac{\eta(2\tau)^8\eta(\tau/2)^8}{\eta(\tau)^8},\label{Cond:Theta8}\\
n=16:\ \Theta_L(q)&=&\frac{\eta(\tau)^{32}}{\eta(2\tau)^{16}}+2^{\frac{k}{2}}\frac{\eta(\tau)^{32}}{\eta(\tau/2)^{16}} +2^{\frac{k}{2}}\frac{\eta(2\tau)^{16}\eta(\tau/2)^{16}}{\eta(\tau)^{16}}\label{Cond:Theta16},
\end{eqnarray}
where $q=e^{2\pi \sqrt{-1}\tau}$ and $\eta(\tau)=q^{1/24}\Pi_{i=1}^\infty(1-q^i)$, which is called the Dedekind $\eta$-function.
The following formulas are well-known:
\eqn
\Theta_L\biggl(\frac{-1}{\tau}\biggr)&=&\biggl(\frac{\tau}{\sqrt{-1}}\biggr)^{n/2}\frac{1}{\sqrt{|L^\circ/L|}}\Theta_{L^\circ}(\tau),\\
\eta\biggl(\frac{-1}{\tau}\biggr)&=&\biggl(\frac{\tau}{\sqrt{-1}}\biggr)^{1/2}\eta(\tau).
\eeqn
Thus, we obtain the theta series of $L^\circ$:
\eqa
n=8:\ \Theta_{L^\circ}(q)&=&\frac{\eta(\tau)^{16}}{\eta(2\tau)^8}+2^{(16-k)/2}\frac{\eta(\tau)^{16}}{\eta(\tau/2)^8} -2^4\frac{\eta(2\tau)^8\eta(\tau/2)^8}{\eta(\tau)^8}\label{Eq:Thetadual8},\\ {}
n=16:\ \Theta_{L^\circ}(q)&=&\frac{\eta(\tau)^{32}}{\eta(2\tau)^{16}}+2^{(32-k)/2}\frac{\eta(\tau)^{32}}{\eta(\tau/2)^{16}} +2^{8}\frac{\eta(2\tau)^{16}\eta(\tau/2)^{16}}{\eta(\tau)^{16}}\label{Eq:Thetadual16}.
\eeqa
In particular, $\Theta_{L^\circ}(q)\in\Z[q^{1/2}]$, which implies that $\langle v,v\rangle\in\Z$ for $v\in L^\circ$.
Therefore $L$ is $2$-elementary totally even.

Suppose that $L$ does not have roots.
Then the coefficients of $q$ in (\ref{Cond:Theta8}) and (\ref{Cond:Theta16}) are zero.
So we have $k=0$ if $n=8$, and $k=8$ if $n=16$.
Hence the determinant of $L$ is $2^8$.
By Proposition \ref{PEX2ETE}, $L$ is isomorphic to either $\sqrt2E_8$ or $\BW$.
\qe
We will indeed check that $\sqrt2E_8$ and $\BW$ satisfy the condition (I) later (see Proposition \ref{PAE8} and Remark \ref{orbiofBW}).

In the rest of this section, we consider the orbit $Q_L$ for an even lattice $L$ of rank $n$ satisfying the following condition:
\begin{enumerate}[(II)]
\item $L$ has no roots and is isomorphic to neither $\sqrt2E_8$ nor $\BW$.\label{Assumption2}
\end{enumerate}
Then, by Proposition \ref{RAM} (2) and (4), $Q_L\subseteq\Orbit$.
Hence, by Lemma \ref{orbit} and the condition (II), we have
\begin{eqnarray}
Q_L\subseteq \{[0]^-,[\lambda]^\pm|\ \lambda\in L^\circ\cap (L/2),\ |(\lambda+L)_2|=2n\}.\label{InEq:Q}
\end{eqnarray}
We show that this inclusion becomes equality.
The extra automorphisms given in Section \ref{SEAVL+} are crucial in the proof of the following theorem.

\bt{lattice3} Let $L$ be an even lattice of rank $n$ satisfying (II).
\begin{enumerate}
\item $Q_L=\{[0]^-,[\lambda]^\pm|\ \lambda\in L^\circ\cap (L/2),\ |(\lambda+L)_2|=2n\}$.
\item If $|Q_L|>1$ then $L$ is obtained by Construction B.
\end{enumerate}
\et

\proof\ Set $Q=\{[0]^-,[\lambda]^\pm|\ \lambda\in L^\circ\cap (L/2),\ |(\lambda+L)_2|=2n\}$.
By (\ref{InEq:Q}), we have $Q_L\subseteq Q$.
If $Q=\{[0]^-\}$ then the assertion is trivial.
So we assume $|Q|>1$.
Let $[\lambda]^+$ be an element of $Q$.
By Proposition \ref{lattice2}, the lattice $L$ is obtained by Construction B associated with an orthogonal basis consisting of vectors in $(\lambda+L)_2$.
By Proposition \ref{extra}, there exists $\sigma\in\Aut(V_L^+)$ such that $[0]^{-,\sigma}=[\lambda]^+$, hence $[\lambda]^+\in Q_L$.
By Proposition \ref{actioncoset}, there exists $f\in\Hom(L,\Z_2)\subset \Aut(V_L^+)$ such that $[\lambda]^{-,f}=[\lambda]^+$.
Therefore $Q\subseteq Q_L$, and $Q=Q_L$.
The arguments above show (2) as well.
\qe

As an application of Theorem \ref{lattice3}, we obtain the following proposition.

\bp{PGOVL+} Let $L$ be an even lattice without roots.
\begin{enumerate}
\item $\Aut(V_L^+)$ is generated by $O(\hat{L})/\langle\theta_{V_L}\rangle$ and the extra automorphisms given in Section \ref{SEAVL+}.
\item $O(\hat{L})/\langle\theta_{V_L}\rangle\subsetneq\Aut(V_L^+)$ if and only if $L$ is obtained by Construction B.
\end{enumerate}
\ep
\proof\ First let us show (1).
The stabilizer $H_L$ of $[0]^-$ is isomorphic to $O(\hat{L})/\langle\theta_{V_L}\rangle$ by Proposition \ref{co1}.
Suppose that $L$ is isomorphic to neither $\sqrt2E_8$ nor $\BW$.
Let $G$ be the subgroup of $\Aut(V_L^+)$ generated by $H_L$ and the extra automorphisms given in Section \ref{SEAVL+}.
By the arguments in the proof of Theorem \ref{lattice3}, $G$ acts transitively on $Q_L$.
Hence we obtain $G=\Aut(V_L^+)$.
We will later check that this assertion holds when $L\cong\sqrt2E_8$ and $\BW$ (see Remark \ref{RgenE8} and \ref{RgenBW}).

Next we prove (2).
If $L$ is realized by Construction B, then $V_L^+$ has the extra automorphisms given by (\ref{extra}), hence $O(\hat{L})/\langle\theta_{V_L}\rangle \subsetneq \Aut(V_L^+)$. 
Conversely if $O(\hat{L})/\langle\theta_{V_L}\rangle\subsetneq\Aut(V_L^+)$ then $|Q|>1$ since the stabilizer of $[0]^-$ is $O(\hat{L})/\langle\theta_{V_L}\rangle$ (cf. Proposition \ref{co1} and \ref{CMain}).
By Theorem \ref{lattice3} (2), $L$ is obtained by Construction B.
\qe

In theorem \ref{lattice3}, we determined $Q_L$.
So we obtain a natural group homomorphism from $\Aut(V_L^+)$ to $Sym(Q_L)$, where $Sym(Q_L)$ is the group of permutations on $Q_L$.
We denote this homomorphism by $\zeta_{V_L^+}$
The shape of $\Aut(V_L^+)$ can be described by the kernel and image of $\zeta_{V_L^+}$.
Since $[0]^-$ belongs to $Q_L$, the kernel $\Ker\zeta_{V_L^+}$ of $\zeta_{V_L^+}$ is a subgroup of the stabilizer of $[0]^-$.
By Theorem \ref{PAVL}, Proposition \ref{co1} and \ref{CMain}, we can determine $\Ker\zeta_{V_L^+}$ in principle.
However, $Sym(Q_L)$ is too large to describe the image $\Imm\zeta_{V_L^+}$ of $\zeta_{V_L^+}$.
So, we consider a certain structure on the orbit $Q_L$ preserved by the action of $\Aut(V_L^+)$.

\bp{lattice5} Let $L$ be an even lattice of rank $n$ satisfying (II).
Set $P=Q_L\cup\{[0]^+\}$.
Then for any $W^1,W^2\in P$, there exists a unique element $W^3$ of $P$ such that $W^1\times W^2=W^3$.
Moreover $P$ forms a vector space over $\F_2$ under this binary operation.
\ep
\proof\ If $Q_L=\{[0]^-\}$ then the assertion is trivial.
We assume $|Q_L|\ge2$.
If $W^1=[0]^+$ or $W^2=[0]^+$ then the first assertion is trivial.
Hence we may assume that $W^1,W^2\in Q_L$.
Then there exists $g\in\Aut(V_L^+)$ such that $W^{1,g}=[0]^-$.
By Theorem \ref{lattice3}, the set $Q_L$ consists of isomorphism classes of irreducible modules of untwisted type.
So $W^{2,g}=[\lambda]^\varepsilon$ for some $\lambda\in L^\circ\cap (L/2)$ and $\varepsilon\in\{\pm\}$.
By Proposition \ref{RFusionFo} (1), the fusion rules $[0]^-\times[\lambda]^\pm=[\lambda]^\mp$ hold for all $\lambda\in L^\circ\cap (L/2)$.
Thus we have $W^1\times W^2=g^{-1}(W^{1,g}\times W^{2,g})=g^{-1}([0]^-\times [\lambda]^\varepsilon)=g^{-1}[\lambda]^{\varepsilon_0}\in P$, where $\{\varepsilon,\varepsilon_0\}=\{\pm\}$.

By Proposition \ref{RFusionFo} (3), the binary operation $\times$ is associative and commutative, hence this operation gives the structure of an abelian group on $P$.
Since $\Aut(V_L^+)$ acts transitively on $Q_L$ and $[0]^-\times[0]^-=[0]^+$, the order of any element of $Q_L$ is $2$.
Therefore $P$ is an elementary abelian $2$-group, namely $P$ forms a vector space over $\F_2$ under the operation $\times$.\qe

By Proposition \ref{RAM} (4), $\Aut(V_L^+)$ preserves the group structure on $P$.
We regard $\Imm\zeta_{V_L^+}$ as a subgroup of the general linear group $GL(P)$ on the vector space $P$.
Set
\begin{eqnarray}
U_L=\Big\{\lambda+L\in (L^\circ\cap (L/2))/L\Big|\ |(\lambda+L)_2|=2n\Big\}\cup\{L\}.\label{Def:U_L}
\end{eqnarray}
Let $\lambda_1+L,\lambda_2+L$ be elements of $U_L$.
By Proposition \ref{RFusionFo} (1), we have $[\lambda_1]^+\times[\lambda_2]^+=[\lambda_1+\lambda_2]^\varepsilon$ for some $\varepsilon\in\{\pm\}$.
It follows form Proposition \ref{lattice5} that $\lambda_1+\lambda_2\in U_L$.
Hence the set $U_L$ is a subspace of $(L^\circ\cap (L/2))/L$.
Since $\AutO(L)=\AutO(L^\circ)$, the group $\AutO(L)$ acts on $U_L$.
Let $\rho_L$ denote the natural group homomorphism from $\AutO(L)$ to $GL(U_L)$.
The image of $\zeta_{V_L^+}$ is described in the following proposition.

\bp{Plattice6} Let $L$ be an even lattice satisfying (II) and let $H_L$ be the stabilizer of the isomorphism class $[0]^-$ of $V_L^-$.
Then $\zeta_{V_L^+}(H_L)$ is of shape $2^{\dim U_L}.\rho_L(\AutO(L))$.
Moreover if $\rho_L$ is surjective then $\zeta_{V_L^+}$ is also surjective.
\ep

\proof\ The first assertion follows from Proposition \ref{co1} and \ref{actioncoset}.
Suppose that $\rho_L$ is surjective.
If $\dim U_L=0$ then $\dim P=1$ and the second assertion is trivial.
So we assume $\dim U_L\ge1$.
Clearly $\zeta_{V_L^+}(H_L)$ is a subgroup of the stabilizer $G$ of $[0]^-$ under the action of $\Imm\zeta_{V_L^+}$ on $P$.
By the first assertion, $\zeta_{V_L^+}(H_L)$ is of shape $2^{\dim U_L}.GL(U_L)$.
The one-point stabilizer under the natural action of $GL_{m+1}(2)$ on non-zero vectors in $\F_2^{m+1}$ is maximal and of shape $2^{m}.GL_m(2)$.
Hence $\zeta_{V_L^+}(H_L)=G$.
By Theorem \ref{lattice3}, Proposition \ref{PGOVL+} (2) and the assumption that $\dim U_L\ge1$, there exists $\sigma\in\Aut(V_L^+)$ such that $\zeta_{V_L^+}(\sigma)\notin G$.
Hence we have $G\subsetneq\Imm\zeta_{V_L^+}\subseteq GL(P)$.
Therefore $\zeta_{V_L^+}$ is surjective.\qe

\subsection{Main results}\label{SMAIN}
We will summarize the results given in Section \ref{SVL+2ETE} and \ref{SAGVL+}.

Let $L$ be an even lattice of rank $n$ and let $Q_L$ be the orbit of the isomorphism class $[0]^-$ of the irreducible $V_L^+$-module $V_L^-$ under the action of $\Aut(V_L^+)$ on the set $S_L$ of all isomorphism classes of irreducible $V_L^+$-modules.

\bt{TMAIN} $(${\bf Main result 1}$)$
\begin{enumerate}
\item The stabilizer $H_L$ of $[0]^-$ is isomorphic to $C_{\Aut(V_L)}(\theta_{V_L})/\langle\theta_{V_L}\rangle$.
\item If $L$ has no roots then $C_{\Aut(V_L)}(\theta_{V_L})/\langle\theta_{V_L}\rangle=O(\hat{L})/\langle\theta_{V_L}\rangle\cong \Hom(L,\Z_2).(\AutO(L)/\langle-1\rangle)$.
\end{enumerate}
\et

\bt{TMAIN0} $(${\bf Main result 2}$)$ Suppose that $L$ is a $2$-elementary totally even lattice with determinant $2^m$.
Then $S_L$ forms an $(m+2)$-dimensional vector space over $\F_2$ under the fusion rules.
Moreover, if $n\in8\Z$ then $S_L$ has an $\Aut(V_L^+)$-invariant quadratic form associated with a non-singular symplectic form.
\et

\bt{TMAIN1} $(${\bf Main result 3}$)$ Suppose that $Q_L$ contains isomorphism classes of irreducible $V_L^+$-modules of twisted type.
Then $L$ is $2$-elementary totally even.
Moreover if $L$ has no roots then it is isomorphic to either $\sqrt2E_8$ or the Barnes-Wall lattice $\BW$ of rank $16$.
\et

\bt{TMAIN2} $(${\bf Main result 4}$)$ Suppose that $L$ does not have roots and that $L$ is isomorphic to neither $\sqrt2E_8$ nor $\BW$.
\begin{enumerate}
\item $Q_L=\{[0]^-,[\lambda]^\pm|\ \lambda\in L^\circ\cap (L/2),\ |(\lambda+L)_2|=2n\}$
\item The set $P=Q_L\cup\{[0]^+\}$ forms a vector space over $\F_2$ under the fusion rules.
\end{enumerate}
\et
See Section \ref{SVL+2ETE} for Theorem \ref{TMAIN0} and Theorem \ref{TMAIN1}, and see Section \ref{SAGVL+} for Theorem \ref{TMAIN} and \ref{TMAIN2}.

The shape of $\Aut(V_L^+)$ for an even lattice $L$ without roots is described in the following way.
Suppose that $L$ is isomorphic to neither $\sqrt2E_8$ nor $\BW$.
By Theorem \ref{TMAIN2} (1), we obtain the orbit $Q_L$ of $[0]^-$.
Then by Theorem \ref{TMAIN2} (2), we have a natural group homomorphism $\zeta_{V_L^+}$ from $\Aut(V_L^+)$ to $GL(P)$.
Clearly, its kernel is a subgroup of the stabilizer $H_L$ of $[0]^-$.
Hence by Theorem \ref{PAVL} and \ref{TMAIN}, we can determine the kernel of $\zeta_{V_L^+}$ in principle.
The image $\Imm\zeta_{V_L^+}$ of $\zeta_{V_L^+}$ is a subgroup of $GL(P)$ and the index of $\zeta_{V_L^+}(H_L)$ in $\Imm\zeta_{V_L^+}$ is equal to $|Q_L|$.
By Proposition \ref{Plattice6} and Theorem \ref{TMAIN}, we determine $\zeta_{V_L^+}(H_L)$.
By group theoretical facts on general linear groups, there are a few possibilities of the shape of $\Imm\zeta_{V_L^+}$.
Thus we can determine the image of $\zeta_{V_L^+}$, and obtain the shape of $\Aut(V_L^+)$.

When $L$ is isomorphic to $\sqrt2E_8$ and $\BW$, the shape of $\Aut(V_L^+)$ will be determined by using Theorem \ref{TMAIN0} in Section \ref{SsSqrt2R} and \ref{SsBW}.
\section{Examples}\label{SExample}
In this section, as an application of the results of the previous section, we determine $\Aut(V_L^+)$ for many important lattices $L$.

\subsection{Even unimodular lattices without roots}\label{SsEUL}
Let $L$ be a positive-definite even unimodular lattice of rank $n$ without roots.
We note that if $n\le16$ then $L$ has roots (cf. \cite{CS}).
Hence $n\ge24$.
By Theorem \ref{TMAIN1} (2), the orbit of $[0]^-$ consists of only $[0]^-$.
Hence the following theorem follows from Theorem \ref{TMAIN}.

\bt{autoofunimodular} Let $L$ be an even unimodular lattice without roots.
Then 
\begin{eqnarray*}
\Aut(V_L^+)\cong\AutO(\hat{L})/\langle\theta_{V_L}\rangle\cong \Hom(L,\Z_2).(\AutO(L)/\langle-1\rangle).
\end{eqnarray*}
\et

\subsection{$\sqrt2R$ for an irreducible root lattice $R$ of type ADE}\label{SsSqrt2R}
Let $R$ be an irreducible root lattice of type $ADE$.
We note that the automorphism group of $R$ was described in Table \ref{TaAR} and $\AutO(R)=\AutO(\sqrt2R)$.
Set $L=\sqrt2R$.
Then $L$ has no roots.
Let $Q_L$ be the orbit of $[0]^-$.
The following lemma is needed to describe $\Aut(V_L^+)$.
\bl{root} 
\begin{enumerate}
\item 
\eqn
|Q_L|=
\left\{\begin{array}{cl}
 \mbox{$1$} & \mbox{${\rm if}\ R=A_n\ (n\neq3),E_6,E_7$},\\
 \mbox{$3$} & \mbox{${\rm if}\ R=A_3,D_n\ (n>4)$},\\
 \mbox{$7$} & \mbox{${\rm if}\ R=D_4$}.
\end{array}
\right.
\eeqn
\item Suppose that $|Q_L|>1$ and that $R\neq E_8$.
Let $U_L$ be the subspace of $(L^\circ\cap(L/2))/L$ given in (\ref{Def:U_L}).
Then the natural group homomorphism $\rho_L:\AutO(L)\to GL(U_L)$ is surjective.
\end{enumerate}
\el
\proof\ The assertion (1) follows from Theorem \ref{lattice3}.
Clearly $|U_L|=(|Q_L|+1)/2$.
If $R\neq D_4$ then it follows from $\dim U_L\le1$ that $\rho_L$ is surjective.
The assertion (2) for the case $R=D_4$ is an easy exercise.
\qe

The shape of $\Aut(V_{\sqrt2R}^+)$ except $R=E_8$ is described in the following theorem.

\bt{autoofroot} Let $R$ be an irreducible root lattice {not of type} $E_8$ and set $L=\sqrt2R$.
Then the shape of $\Aut(V_L^+)$ is described as follows:
\begin{enumerate}[{\rm (i)}]
\item If $R=A_n$ $(n\neq3)$, $E_6$ or $E_7$ then $\Aut(V_L^+)=\AutO(\hat{L})/\langle\theta_{V_L}\rangle$.
\item If $R=A_3$ then $\Aut(V_L^+)$ is of shape $(2^2:Sym_4).Sym_3$.
\item If $R=D_4$ then $\Aut(V_L^+)$ is of shape $(2^4:Sym_4).GL_3(2)$.
\item If $R=D_n$ $(n\ge5)$ then $\Aut(V_L^+)$ is of shape $(2^{2n-3}:Sym_n).Sym_3$.
\end{enumerate}
\et
\proof\ If $|Q_L|=1$ then $\Aut(V_L^+)\cong \AutO(\hat{L})/\langle\th_{V_L}\rangle$ by Theorem \ref{TMAIN}.
Hence (i) follows from Lemma \ref{root} (1).

Suppose $|Q_L|>1$.
By Theorem \ref{TMAIN2} (2), the set $P=Q_L\cup\{[0]^+\}$ forms a vector space over $\F_2$ under the fusion rules.
By Proposition \ref{Plattice6} and Lemma \ref{root} (2), the group homomorphism $\zeta_{V_L^+}:\Aut(V_L^+)\to GL(P)$ is surjective.
Let us determine its kernel.
By Remark \ref{Rsplit}, $O(\hat{L})\cong\Hom(L,\Z_2):O(L)$.
Since $[0]^-$ belongs to $Q_L$, $\Ker\zeta_{V_L^+}$ is a subgroup of the stabilizer of $[0]^-$.
By Theorem \ref{TMAIN}, 
\begin{eqnarray}
\Ker\zeta_{V_L^+}\subset\AutO(\hat{L})/\langle\theta_{V_L}\rangle\cong \Hom(L,\Z_2):(\AutO(L)/\langle-1\rangle).\label{Eq:ShapeR1}
\end{eqnarray}
By Proposition \ref{Plattice6},
\begin{eqnarray}
\zeta_{V_L^+}(\AutO(\hat{L})/\langle\theta_{V_L}\rangle))\cong 2^{\dim P-1}:GL_{\dim P-1}(2)\label{Eq:ShapeR2}.
\end{eqnarray}
By comparing (\ref{Eq:ShapeR1}) and (\ref{Eq:ShapeR2}), we obtain the shape of $\Ker\zeta_{V_L^+}$.
Hence (ii), (iii) and (iv) follow.
\qe

\bn{D4}{\rm The automorphism group of $V_{\sqrt2D_4}^+$ was determined in \cite{MM}.
In particular, $\Aut(V_{\sqrt2D_4}^+)\cong 2^6:(GL_3(2)\times Sym_3)$.}
\en

The group $\Aut(V_{\sqrt2E_8}^+)$ was already obtained in \cite{Gr2} (see \cite{KM} for another proof) and the classification of $3$-transposition groups is crucial in the papers.
We will determine the automorphism group of $V_{\sqrt2E_8}^+$ without using the classification.

We recall the properties of the lattice $N=\sqrt2E_8$.
By Example \ref{REXCB} (2), $N$ is obtained by Construction B.
By Proposition \ref{PEX2ETE} (1), $N$ is a $2$-elementary totally even lattice of rank $8$ with determinant $2^8$.
In particular, $N^\circ/N\cong 2^8$.
It is easy to see that the type of the quadratic form $N^\circ/N\to\F_2$, $\lambda+N\mapsto\langle\lambda,\lambda\rangle\mod2$ is plus (cf. \cite{CS}).
By Table \ref{TaAR}, $\AutO(N)\cong 2.O^+_8(2)$.
By \ref{exactauto2}, $\Aut(\hat{N})\cong 2^8.2.O^+_8(2)$.

Let $S_N$ be the set of all isomorphism classes of irreducible $V_N^+$-modules.
By Remark \ref{symonlattice} and Theorem \ref{TMAIN0}, $S_N$ forms a $10$-dimensional vector space over $\F_2$ under the fusion rules and has the $\Aut(V_N^+)$-invariant quadratic form $q_N$ of plus type.
Let $O(S_N)$ denote the orthogonal group on $S$ associated with $q_N$ and let $\xi_N:\Aut(V_{N}^+)\to O(S_N)\cong O^+_{10}(2)$ be the natural group homomorphism.

Let us determine both the kernel and image of $\xi_N$.
It is well known that $\AutO(N)/{\B-1\rangle}\cong O^+_8(2)$ acts faithfully on $N^\circ/N$.
Since $N^\circ=N/2$, we have $\{h\in \Hom(N,\Z_2)|\ h(2N^\circ)=0 \}=\{1\}$.
By Proposition \ref{kernelofauto3}, the kernel of $\xi_N$ is trivial.
Hence $\Aut(V_N^+)\subseteq O(S_N)$.
By Proposition \ref{extra}, there is an automorphism $\sigma$ of $V_N^+$ such that $\sigma\notin \AutO(\hat{N})/\langle\th_{V_N}\rangle$.
Hence we have $\langle \AutO(\hat{N})/\langle\theta_{V_N}\rangle,\sigma\rangle\subseteq \Aut(V_N^+)\subseteq O(S_N)$.
Since $\AutO(\hat{N})/\langle\theta_{V_N}\rangle\cong 2^8:O^+_8(2)$ is a maximal subgroup of $O^+_{10}(2)$ (cf.\ \cite{ATLAS}), we have $\Aut(V_N^+)\cong O(S_N)\cong O^+_{10}(2)$.

\bt{autoofE8} {\rm \cite{Gr2,KM}} The automorphism group of $V_{\sqrt2E_8}^+$ is isomorphic to $O^+_{10}(2)$.\et

\br{RgenE8}{\rm The automorphism group of $V_{N}^+$ is generated by $O(\hat{N})/\langle\theta_{V_{N}}\rangle$ and the extra automorphisms given in Section \ref{SEAVL+}, where $N=\sqrt2E_8$.}
\er

Now we study the orbit $Q_N$ of $[0]^-$ more precisely.
According to \cite{ATLAS}, the orbit decomposition of $S_N$ under the action of the orthogonal group $O(S_N)$ is $2^{10}=1+496+527$, where $1$, $496$, $527$ mean the zero vector, the set of $496$ non-isotropic vectors, the set of $527$ non-zero isotropic vectors respectively.
By the definition of the quadratic form $q_N$ given in Section \ref{SVL+2ETE}, $[0]^-$ is non-zero isotropic.
Hence we obtain the following proposition.

\bp{PAE8} The orbit $Q_{N}$ is the set of all non-zero isotropic vectors in $(S_N,q_N)$, namely
\begin{eqnarray*}
Q_{N}=\{[0]^-,[\lambda]^\pm,\ [\chi]^-|\ \lambda\in N^\circ\cap (N/2),\ |(\lambda+N)_2|=16,\ \chi\in X(\Cent(\hat{N}/K_N))\},
\end{eqnarray*}
where $N=\sqrt2E_8$.
In particular, the cardinality of $Q_{N}$ is $527$.
\ep

The results of this subsection are summarized in Table \ref{autorootlatticeVOA}.

\begin{table}
\caption{Automorphism group of $V_{\sqrt2R}^+$}\label{autorootlatticeVOA}
\begin{center}
\begin{tabular}{|c|| c|c|} 
\hline 
$R$&$\Aut(V_{\sqrt2R}^+)$& $|Q_{\sqrt2R}|$\\ \hline
$A_1$& $\Z_2$&$1$\\ \hline
$A_n$ $(n\neq1,3)$&$2^n:Sym_{n+1}$&$1$\\ \hline
$A_3$&$(2^{2}:Sym_4).Sym_3$&$3$\\ \hline
$D_4$&$(2^{4}:Sym_4).GL_3(2)$&$7$\\ \hline
$D_n$ $(n\ge5)$&$(2^{2n-2}:Sym_n).Sym_3$&$3$\\ \hline
$E_6$&$2^6:U_4(2):2$&$1$\\ \hline
$E_7$&$2^7:Sp_6(2)$&$1$\\ \hline
$E_8$&$O_{10}^+(2)$&$527$\\
\hline
\end{tabular}
\end{center}
\end{table}

\subsection{Barnes-Wall lattice of rank $16$}\label{SsBW}
We recall the properties of the Barnes-Wall lattice $\BW$ of rank $16$.
By Example \ref{REXCB} (3), $\BW$ is obtained by Construction B.
By Proposition \ref{PEX2ETE} (2), $\BW$ is a $2$-elementary totally even lattice of rank $16$ with determinant $2^8$ without roots.
By Table \ref{TaAR}, we have $\AutO(\BW)\cong2^{1+8}_+\cdot \Omega^+_8(2)$.
It is easy to see that the type of the quadratic form $\BW^\circ/\BW\to\F_2$, $\lambda+\BW\mapsto\langle\lambda,\lambda\rangle\mod2$ is plus (cf. \cite{CS}).

Let $S_{\BW}$ be the set of all isomorphism classes of irreducible $V_{\BW}^+$-modules.
By Remark \ref{symonlattice} and Theorem \ref{TMAIN0}, $S_{\BW}$ forms a $10$-dimensional vector space over $\F_2$ under the fusion rules and has the $\Aut(V_{\BW}^+)$-invariant quadratic form $q_{\BW}$ of plus type.
Let $O(S_{\BW})$ denote the orthogonal group on $S$ associated with $q_{\BW}$ and let $\xi_{\BW}:\Aut(V_{\BW}^+)\to O(S_{\BW})\cong O^+_{10}(2)$ be the natural group homomorphism.

First we determine the kernel $F$ of $\xi_{\BW}$.
Let $H_{\BW}$ be the stabilizer of the isomorphism class $[0]^-$ of $V_{\BW}^-$.
Then $F$ is a normal subgroup of $H_{\BW}$.
By Theorem \ref{TMAIN}, we have
\begin{eqnarray}
H_{\BW}=C_{\Aut(V_{\BW})}(\theta_{V_{\BW}})/\langle\theta_{V_{\BW}}\rangle\cong\Hom(\BW,\Z_2).(\AutO(\BW)/\langle-1\rangle).\label{Eq:ShapeH}
\end{eqnarray}
In particular, $H_{\BW}$ is of shape $2^{16}.2^8.\Omega_{8}^+(2)$.
We set $F_1=F\cap\Hom(\BW,\Z_2)$.
By Proposition \ref{kernelofauto3}, $F_1\cong 2^8$ and $F/F_1\cong 2^8$.
In order to determine the structure of $F$, we use the automorphism $\sigma\in\Aut(V_{\BW}^+)$ given in Section \ref{SEAVL+}.
It is easy to see that $\sigma^{-1} F_1\sigma \neq F_1$ and the group $F_1$ is a subgroup of the center $\Cent(F)$ of $F$. 
Since $\AutO(\hat{\Lambda}_{16})$ acts irreducibly on $F/F_1$ and $\sigma$ does not normalize $F_1$, we have $F=\Cent(F)$.
Since $F_1$ is an elementary abelian $2$-group, $F$ is generated by elements of order $2$.
Hence we obtain $F\cong2^{16}$.

Next we determine the image of $\xi_{\BW}$.
Set $G=H_{\BW}/F$.
By (\ref{Eq:ShapeH}) and the shape of $F$, the shape of $G$ is $2^8.\Omega_8^+(2)$.
Since $\sigma\notin H_{\BW}$, we have $\langle G,\sigma\rangle\subset O(S_{\BW})$.
Since $2^8.\Omega_8^+(2)$ is a maximal subgroup of $\Omega_{10}^+(2)$ and $|O^+_{10}(2):\Omega_{10}^+(2)|=2$, the image of $\xi_{\BW}$ is of shape $\Omega_{10}^+(2)$ or $O^+_{10}(2)$.

Now we consider the action of $O(S_{\BW})$ on $S_{\BW}$.
Similar arguments as in the previous section show that the cardinality of the orbit $Q_{\BW}$ of $[0]^-$ is less than or equal to $527$, which is the number of non-zero isotropic vectors in $S_{\BW}$.
Hence we have
\begin{eqnarray*}
|\Omega^+_{10}(2)|\le |\Imm\xi_{\BW}|\le 527\cdot|G|.
\end{eqnarray*}
Since $527\cdot|G|=|\Omega_{10}^+(2)|$, we have $\Imm\xi_{\BW}=\langle G,\sigma\rangle\cong \Omega_{10}^+(2)$.
In particular, the cardinality of $Q_{\BW}$ is $527$.
Therefore we obtain {an} exact sequence
\begin{eqnarray}
1\to F\to \Aut(V_{\BW}^+)\to\langle G,\sigma\rangle\to 1.\label{Ex:BW}
\end{eqnarray}

Finally we show that this exact sequence is non-split.
According to Theorem 1 in $\cite{Gr1}$, the shape of the group $B=\AutO(\BW)/\langle-1\rangle$ is $2^8.\Omega^+_8(2)$ and it is non-split.
{Since $\Omega^+_8(2)$ is simple}, $B$ does not have a subgroup of shape $\Omega^+_8(2)$.
Let $A$ be the normal subgroup of $B$ of shape $2^8$.
Set $K=\Hom(\BW,\Z_2)\subset H_{\BW}$.
Then $H_{\BW}/K\cong B$.
Suppose that (\ref{Ex:BW}) is split.
Let $C$ be a complement to $F$ in $\Aut(V_{\BW}^+)$.
It is easy to see that $K(C\cap H_{\BW})/K\cong B/A$.
This implies that $H_{\BW}/K$ has a subgroup of the shape $\Omega^+_8(2)$, which is a contradiction.
Hence the exact sequence (\ref{Ex:BW}) is non-split.

Summarizing the results above, we have the following theorem.
\bt{nonsplitofG} The shape of the automorphism group of $V_{\BW}^+$ is $2^{16}\cdot \Omega_{10}^+(2)$.
\et

\br{RgenBW}{\rm The automorphism group of $V_{\BW}^+$ is generated by $O(\hat{\Lambda}_{16})/\langle\theta_{V_{\BW}}\rangle$ and the extra automorphisms given in Section \ref{SEAVL+}.}
\er

\br{orbiofBW}{\rm As in Proposition \ref{PAE8}, the orbit $Q_{\BW}$ of $[0]^-$ is the set of all non-zero isotropic vectors in $(S_{\BW},q_{\BW})$, namely
\begin{eqnarray*}
Q_{\BW}=\{[0]^-,[\lambda]^\pm,\ [\chi]^+|\ \lambda\in \BW^\circ\cap (\BW/2),\ |(\lambda+\BW)_2|=32,\ \chi\in X(\Cent(\hat{\Lambda}_{16}/K_{\BW}))\}.
\end{eqnarray*}
In particular the cardinality of $Q_{\BW}$ is 527.}
\er

\end{document}